\documentclass[pdflatex,sn-mathphys-num]{sn-jnl}

\usepackage{graphicx}
\usepackage{multirow}
\usepackage{amsmath,amssymb,amsfonts}
\usepackage{amsthm}
\usepackage{mathrsfs}
\usepackage{xcolor}
\usepackage{booktabs}
\usepackage{algorithm}
\usepackage{algorithmicx}
\usepackage{algpseudocode}

\theoremstyle{thmstyleone}
\newtheorem{thm}{Theorem}
\newtheorem{theorem}[thm]{Theorem}
\newtheorem{lem}[thm]{Lemma}
\newtheorem{proposition}[thm]{Proposition}
\newtheorem{corollary}[thm]{Corollary}
\theoremstyle{thmstyletwo}

\theoremstyle{thmstylethree}
\newtheorem{definition}{Definition}
\DeclareMathOperator*{\Limsup}{Lim\,sup}
\begin{document}
\title[Exact penalty and adaptive discretization for nonsmooth SIP]{Exact penalty methods with adaptive discretization for nonsmooth semi-infinite optimization}

\author[1]{\fnm{Giancarlo} \sur{Bigi}}\email{giancarlo.bigi@unipi.it}
\author*[2]{\fnm{Riccardo} \sur{Tomassini}}\email{riccardo.tomassini@uniroma1.it}

\affil[1]{\orgdiv{Department of Computer Science}, \orgname{University of Pisa}, \orgaddress{\city{Pisa}, \country{Italy}}}
\affil*[2]{\orgdiv{Department of Computer, Control and Management Engineering}, \orgname{Sapienza University of Rome}, \orgaddress{\city{Rome}, \country{Italy}}}

\abstract{The paper considers a nonsmooth semi-infinite optimization problem in which the constraints are reformulated through a value function. Stationary points of the cutting surface approximations are studied, showing that stationarity is preserved in the limit under boundedness of the associated sequence of multipliers. Building on these results, an algorithm that combines the cutting surface approximations with a penalty scheme is devised, together with some proximal  variants. An exact penalty result is established under suitable global conditions, and numerical results are  reported to compare the  methods with existing exact penalty and cutting surface algorithms.}


\maketitle

\section{Introduction}
Semi-infinite programs are optimization problems with finitely many decision variables and infinitely many constraints. This structure naturally arises when a constraint must be satisfied for every value of a parameter varying over a continuous domain. As a result, semi-infinite models appear in robotics, engineering design, economic planning,  optimal control, approximation theory, and robust optimization (see, e.g. \cite{reemtsen_semi-infinite_1998,lopez_semi-infinite_2007,goberna_robust_2014} and the references therein).
 In this work, we consider the semi-infinite optimization problem
\begin{align*}
\min\quad &
 f(x)\\
\notag\
\textnormal{s.t.} \quad &
 g(x,y)\leq 0
\qquad
\forall y\in Y,\\
& x \in X,
\tag{SIP}\label{sip}
\end{align*}
where $X\subseteq\mathbb{R}^{n}$ is  closed and convex, and
$Y\subseteq\mathbb{R}^{m}$ is  compact. The function
$f:\mathbb{R}^{n}\to\mathbb{R}$ is locally Lipschitz continuous, while
$g:\mathbb{R}^{n}\times\mathbb{R}^{m}\to\mathbb{R}$ is continuous and
$g(\cdot,y)$ is locally Lipschitz continuous for every $y\in Y$.

Semi-infinite programs have been extensively studied from both theoretical and computational perspectives, leading to several approximation and reformulation techniques for handling the infinite constraints (see, e.g.
\cite{reemtsen_semi-infinite_1998,lopez_semi-infinite_2007,hettich_httpwwwjstororg_1993,shapiro_semi-infinite_2009,djelassi_recent_2021}).
A major class of numerical methods is based on discretization, where only a finite subset of the constraints is considered. Grid-based methods use increasingly dense finite subsets
(see, e.g. \cite{still_discretization_2001,turan_bounding-focused_2025}), while cutting surface methods use adaptive discretizations that refine the approximation of the feasible region by adding or removing selected constraints according to the current iterate
(see, e.g. \cite{blankenship_infinitely_1976,tsoukalas_feasible_2011,mehrotra_cutting_2014,CerulliOustryDAmbrosioLiberti2022}).
Another class of approaches reformulates the semi-infinite problem into a standard optimization problem  through an integral function that measures the overall violation of the  constraints. Integral reformulations have been extensively studied and are often combined with penalty and smoothing techniques
(see, e.g. \cite{conn_exact_1987,liu_convergence_2013,lin_new_2014,xu_solving_2014}).
In the convex setting, Remez-type methods combine penalty techniques with discretization schemes
\cite{auslender_comparative_2015}.
Other approaches for convex semi-infinite problems include primal-dual methods
(see, e.g. \cite{ItoLiuTeo2000,WeiHaskellZhao2020.1})
and bundle methods
(see, e.g. \cite{PangLvWang2016,PangWuWangWu2020}).
For smooth problems, another strategy exploits the KKT conditions of the semi-infinite program  to derive finite systems of equations and inequalities that can be solved with ad hoc numerical methods
(see, e.g. \cite{WuLiQiZhou2005,JinLingShen2015}).
For general nonconvex semi-infinite problems, global optimization methods have also been developed by constructing and iteratively refining valid lower and upper bounds through suitable relaxations of the constraints
(see, e.g. \cite{BhattacharjeeGreenBarton2005,floudas_adaptive_2008,Mitsos2011,shiu_relaxed_2012}).

In this work, the infinite constraints are replaced by the value function, that measures the maximum violation of the constraints  at each point. While  reformulations with a unique constraint have already been studied, they often rely on integral-type constructions to obtain smooth functions. The  value function instead avoids numerical integration over the index set $Y$, therefore it can  better handle problems where $Y$ lies in a  high dimensional space. An approximation of the value function is employed  by measuring the maximum violation of the constraints over  a finite subset of $Y$. The resulting optimization problem is a relaxation of the original semi-infinite problem while being easier to solve. This is the idea behind cutting surface methods, which iteratively solve such approximate problems and  update the  subset of the constraints according to the current iterate. In particular, existing methods often rely  on sequences of global minimizers of the approximate problems to recover an optimal solution of the original semi-infinite problem in the limit. This can be restrictive in  nonsmooth and nonconvex settings, where global minimizers may be difficult to compute. For this reason, our analysis focuses on Clarke stationary points of a generic sequence of approximate problems. The limiting behavior  of the associated approximated value functions and Clarke generalized gradients are studied through the Kuratowski superior limit. This allows proving that stationarity is preserved in the limit whenever the corresponding sequence of Lagrange multipliers is bounded.

To devise numerical methods based on exact penalization, the positive part of the approximate value function is added to the objective and weighted by a penalty parameter. This leads to an algorithm that computes a stationary point of the current penalized approximation at each iteration. The penalty parameter is increased adaptively to promote feasibility, while the cutting surface approximation is refined by adding new constraints. Convergence to stationary points of the semi-infinite problem is established when the penalty sequence remains bounded, or when  a cluster point satisfies a suitable constraint qualification. Then,
 proximal variants are
introduced to improve convergence by adding a quadratic regularization term in the objective. The proximal structure allows the objective function to be locally linearized and inactive constraints to be removed, leading to simpler approximate problems. Since the convergence analysis relies on the boundedness of the  penalty parameters, an exact penalty result is established under compactness of the feasible region $X$ and a global constraint qualification.

The article is organized as follows. Section~2 recalls the main tools from nonsmooth and set-valued analysis that are needed to define stationarity and study the convergence. Section~3 develops the  value function reformulation, introduces its cutting surface approximations and the  penalty relaxations, and analyzes the relations between the  corresponding stationarity conditions. Section~4 presents the penalty and proximal cutting surface algorithms, establishes their convergence properties, and derives the exact penalty result under suitable additional assumptions.  Section~5 reports the  results of some numerical tests, that compare the  algorithms presented in Section~4 with each other and with existing exact penalty and cutting surface methods.

\section{Mathematical toolbox}
\label{sec:toolbox}

This section collects the notions from nonsmooth and set-valued analysis that are needed  throughout the paper. In fact,  Clarke generalized gradients and normal cones are exploited to achieve first-order stationarity conditions, while the convergence analyses  further require notions of continuity (closedness) and limits for set-valued function.

Let $h:\mathbb{R}^p\to\mathbb{R}$ be locally Lipschitz continuous. The Clarke generalized directional derivative of $h$ at $x\in\mathbb{R}^p$ in the direction $d\in\mathbb{R}^p$ is defined by
\begin{align*}
h^{\circ}(x;d)
&:=
\limsup_{\substack{z\to x\ t\downarrow 0}}
\frac{h(z+td)-h(z)}{t},
\end{align*}
while
the standard directional derivative is given by
\begin{align*}
h'(x;d)
&:=
\lim_{t\downarrow 0}
\frac{h(x+td)-h(x)}{t},
\end{align*}
provided that the limit exists. Unlike the Clarke generalized directional derivative, this limit is not guaranteed to exist for a locally Lipschitz function. When it exists and  coincides with  the   Clarke generalized directional derivative for every direction, i.e.,
\begin{align*}
h'(x;d)
&=
h^{\circ}(x;d)
\qquad
\forall d\in\mathbb{R}^p,
\end{align*}
then we say that function $h$ is Clarke regular at $x$.
Clarke regularity ensures that the standard and generalized directional derivatives provide the same  first-order information.
Stationarity conditions are then formulated in terms of both the Clarke generalized gradient and the normal cone. We define the Clarke generalized gradient  of $h$ at $x$
\begin{align*}
\partial h(x)
&:=
\left\{
\xi\in\mathbb{R}^p
:
\xi^\top d\leq h^{\circ}(x;d)
\quad
\text{for all }d\in\mathbb{R}^p
\right\}
\end{align*}
and the normal cone of $X$ at $x\in X$
\begin{align*}
\mathcal{N}_X(x)
&:=
\left\{
v\in\mathbb{R}^p
:
v^\top(z-x)\leq 0
\quad
\text{for all }z\in X
\right\}.
\end{align*}
Both $\partial h(x)$ and $\mathcal{N}_X(x)$ are closed and convex sets, while the former is also compact thanks to the local lipschitzianity of $h$ \cite[Proposition 2.1.1]{clarke_optimization_1990}. Moreover, they can be both seen as set-valued maps.  A notion of continuity is needed for set-valued maps in order to analyze their limits, which is given by outer semi-continuity. In particular, a set-valued map
$F:\mathbb{R}^p\rightrightarrows\mathbb{R}^q$ is outer semi-continuous if it has a closed graph, i.e., the following implication
\begin{align*}
x_k\to x_\infty,\qquad
v_k\to v_\infty,\qquad
v_k\in F(x_k)
\quad\Longrightarrow\quad
v\in F(x_\infty)
\end{align*}
holds.  While additional assumptions are needed to guarantee   that  the map
$\partial h:X\rightrightarrows \mathbb{R}^p$ has a closed graph,  the normal cone map
$\mathcal{N}_X:X\rightrightarrows \mathbb{R}^p$ is always outer semi-continuous if $X$ is convex and closed \cite[Proposition 6.6]{RockafellarWets1998}.
 When  $\partial h$ is outer semi-continuous at $x$, the limit
\begin{align*}
\Limsup_{k\to\infty}\partial h(x_k)
&\subseteq
\partial h(x)
\end{align*}
holds true whenever  $x_k \rightarrow x$,  where the Kuratowski superior limit of a    sequence $\{A_k\}_k$ of  subsets of $\mathbb{R}^p$ is
\begin{align*}
\Limsup_{k\rightarrow \infty}(A_k) := \{z \in \mathbb{R}^p : \liminf_{k \to \infty} d(z, A_k) = 0\}.
\end{align*} Similarly, if $X$ is convex and closed, we get
\begin{align*}
\Limsup_{k\to\infty}\mathcal{N}_X(x_k)
&\subseteq
\mathcal{N}_X(x).
\end{align*}

Finally, let us introduce some operations on sets and their sequences. The Minkowski sum of two sets
$A$ and $B$ is
\begin{align*}
A+B
&:=
\left\{
a+b
:
a\in A,\ b\in B
\right\}.
\end{align*}
Given  $\Lambda\subseteq\mathbb{R}$ and $A\subseteq\mathbb{R}^p$, the set of all possible  scalar multiplication  between elements of the two sets is
\begin{align*}
\Lambda A
&:=
\left\{
\lambda a
:
\lambda\in\Lambda,\ a\in A
\right\}.
\end{align*}
To describe the behavior of Kuratowski limits under set operations, let $\{B_k\}_k$  and $\{\Lambda_k\}_k$ be  sequences of  subsets of respectively $\mathbb{R}^p$ and $\mathbb{R}$ with  nonempty  Kuratowski superior limits. Assume that $\{A_k\}_k$ and $\{\Lambda_k\}_k$ are uniformly bounded.
Then, it holds that
\begin{align*}
\Limsup_{k\to\infty}(A_k+B_k)
&\subseteq
\Limsup_{k\to\infty}A_k
+
\Limsup_{k\to\infty}B_k
\end{align*}
and
\begin{align*}
\Limsup_{k\to\infty}(\Lambda_k A_k)
&\subseteq
\left(\Limsup_{k\to\infty}\Lambda_k\right)
\left(\Limsup_{k\to\infty}A_k\right).
\end{align*}

\section{Stationary points}
{ The infinite family of constraints is replaced by a value function that measures the maximum  violation of the constraints at each point. Since the resulting optimization problem is generally nonsmooth and nonconvex,  stationarity conditions are used as a more suitable target for optimization. The  evaluation of the value function requires taking into account   the whole index set, which could be computationally demanding. For this reason, cutting surface methods are used to approximate it by using only a possibly finite subset of indices. For each resulting
approximate problem, a penalized problem is considered by incorporating the constraint violation into the objective function. Finally, the relations between stationary points of the original problem and  its approximations are studied to establish  convergence results for the methods.}

The following properties are assumed to hold throughout the paper:
\begin{itemize}
\item[\textbf{(A1)}]
$X\subseteq\mathbb{R}^n$ is
closed and convex, and $Y\subseteq\mathbb{R}^m$ is compact.

\item[\textbf{(A2)}]
 $f:\mathbb{R}^n\to\mathbb{R}$ is locally Lipschitz continuous and Clarke regular on  $X$ and the map $\partial f:X\rightrightarrows \mathbb{R}^n$ is outer semi-continuous.

\item[\textbf{(A3)}]
 $g:\mathbb{R}^n\times\mathbb{R}^m\to\mathbb{R}$ is continuous and  $g(\cdot,y)$ is locally Lipschitz continuous and Clarke regular on $X$ for every $y \in Y$. The set-valued map
$\partial_x g:X\times Y\rightrightarrows \mathbb{R}^n$
is outer semi-continuous
and locally bounded.
\end{itemize}
Under assumptions (A1) and (A3) the  function $g(\cdot,y)$ is uniformly Lipschitz on every compact subset of $X$. Indeed, given any compact  $V \subseteq X$ there exists a constant $H$ such that
\begin{align*}
|g(x,y)-g(z,y)|
\leq
H||x-z||
\end{align*}
holds
for all $x,z\in V$ and $y\in Y$.
\subsection{Value function and single level reformulation}

The infinitely many constraints in  \eqref{sip} can be equivalently represented by means of the value function
\begin{align*}
\varphi(x)
&:=
\max_{y\in Y}g(x,y).
\end{align*}
The continuity of $g$ and the compactness of $Y$ ensure that the maximum is attained for every given $x\in X$, so that the value function is well-defined at every point. Moreover, since  $g(\cdot,y)$ is uniformly Lipschitz on every compact subset,
$\varphi$ is also locally Lipschitz continuous and Clarke regular on $X$. Its Clarke generalized gradient is therefore nonempty, compact, and convex at every point.

 The value function allows recasting  \eqref{sip} as
\begin{align}
\min \quad
& f(x)
\notag\ \\
\textnormal{s.t.}\quad
&  \varphi(x)\leq 0, \quad x\in X, \quad
\tag{VP}\label{VPreformulation}
\end{align}
so that  the infinite family of constraints is replaced with a single constraint, which is
generally nonsmooth and nonconvex.
As a consequence, computing global minimizers is generally difficult, and first-order stationary points are therefore a more tractable target for numerical methods. Moreover,  first-order stationarity is a necessary condition for local optimality under a suitable constraint qualification.
Although the value function may be nondifferentiable, stationarity can still be formulated in terms of  generalized derivatives and the normal cone.

\begin{definition}
A point $\bar{x}\in X$ is a stationary point for \eqref{VPreformulation}, or equivalently for \eqref{sip}, if it is feasible   and there exists a multiplier $\lambda\geq 0$ such that the pair $(\bar{x},\lambda)$ satisfies
\begin{align*}
& 0
\in
\partial f(\bar{x})
+
\lambda\partial\varphi(\bar{x})
+
\mathcal{N}_X(\bar{x}),\\
&\lambda\varphi(\bar{x})
=0.
\end{align*}
\end{definition}

A suitable constraint qualification is needed to guarantee that every local minimizer satisfies the above  stationarity conditions. We rely on  the no-nonzero-abnormal-multiplier constraint qualification, which has already been  used in the analysis of value function reformulations of semi-infinite programs (see, e.g. \cite{xu_solving_2014}).

\begin{definition}
A feasible point $\bar{x}\in X$ satisfies the no-nonzero-abnormal-multiplier constraint qualification, shortly NNAMCQ, for \eqref{VPreformulation} if either $\varphi(x)<0$ or $\varphi(x)=0$ and
\begin{align*}
0
\notin
\partial\varphi(\bar{x})
+
\mathcal{N}_X(\bar{x}).
\end{align*}
\end{definition}

The condition automatically holds at points that satisfy the inequality constraint strictly and  at  stationary points whose set of multipliers is bounded. Local minimizers for \eqref{VPreformulation} that satisfy this regularity condition are  stationary points.

\begin{lem}[{{\cite[Corollary 2.1]{xu_solving_2014}}}]

Let $\bar{x}$ be a local minimizer for \eqref{sip}. If NNAMCQ holds at $\bar{x}$, then $\bar{x}$ is a stationary point for  \eqref{sip}.
\end{lem}

The standard NNAMCQ is formulated only at feasible points.
In order to devise penalty schemes, which generate infeasible points, we extend this regularity condition beyond the feasible set.
\begin{definition}
A point $\bar{x}\in X$ satisfies the extended no-nonzero-abnormal-multiplier constraint qualification, shortly  eNNAMCQ, for \eqref{VPreformulation} if either $\varphi(\bar{x})<0$ or $\varphi(\bar{x})\ge0$ and
\begin{align*}
0
\notin
\partial\varphi(\bar{x})
+
\mathcal{N}_X(\bar{x}).
\end{align*}
\end{definition}

Note that  the value function $\varphi$ is also convex  when $g(\cdot,y)$ is convex for every $y\in Y$, and it is easier to check whether the constraint qualification is satisfied. Indeed, in this setting eNNAMCQ holds at every point  whenever a Slater point $x^*$ exists, that is, any point $x^*\in X$ such that $\varphi(x^*)<0$.

\subsection{Cutting surface approximation}

The value-function reformulation replaces the infinite family of  constraints with a unique  inequality constraint, but the evaluation of the value function at a single point still requires solving a  maximization problem over the whole index set $Y$. Cutting surface methods tackle this issue  by considering only a possibly finite subset of constraints.

 Let us define the approximate value function
\begin{align*}
\varphi_U(x)
&:=
\sup_{y\in U}g(x,y),
\end{align*}
where $U\subseteq Y$  ,  so that only the constraints selected  by the index set  $U$ are considered. Therefore,
 the
 corresponding approximate problem
\begin{align*}
\min\quad
& f(x)\\
\tag{VP$_U$} \label{VPU}\
\textnormal{s.t.}\quad
&  \varphi_U(x)\leq 0, \quad  x \in X.
\end{align*}
can be considered.
Clearly, the inequality
$\varphi_U(x)\leq\varphi(x)$ holds at every $x$,
hence every feasible point for   \eqref{sip} is feasible for  \eqref{VPU} as well. Thus, \eqref{VPU} is a relaxation of the original problem as its optimal value is no larger than that of \eqref{sip}.

When $U$ is closed,
  $\varphi_U$ shares the main regularity properties of the original value function, that is,
  $\varphi_U$  is locally Lipschitz continuous and Clarke regular. Therefore, the  Clarke generalized gradient of $\varphi_U$ can  be characterized in terms of appropriate  partial generalized gradients of the constraint $g$ (see \cite[Theorem 2.8.2]{clarke_optimization_1990}).

\begin{proposition}\label{ClarkeGrad}

Let $U\subseteq Y$  nonempty and closed. Then, $\varphi_U$ is Clarke regular at every $x\in X$, and
\begin{align*}
\partial\varphi_U(x)
&=
\operatorname{co}
\left\{
\xi\in\partial_x g(x,y)
\ : \
y\in \arg\max_{y\in U}g(x,y)
\right\},
\end{align*}
where $\operatorname{co}$ denotes the convex hull. Moreover, $\partial\varphi_U(x)$ is compact for any $x\in X$.
\end{proposition}

The regularity properties of the approximate value functions allows defining  stationarity in a similar fashion of   \eqref{VPreformulation}. Therefore,  a feasible point $\bar{x}\in X$ for  \eqref{VPU}  is stationary  if there  exists a multiplier $\lambda_U\geq 0$ such that the pair $(\bar{x},\lambda_U)$ satisfies
\begin{align*}
&0
\in
\partial f(\bar{x})
+
\lambda_U\partial\varphi_U(\bar{x})
+
\mathcal{N}_X(\bar{x}),\\
&\lambda_U\varphi_U(\bar{x})
=0.
\end{align*}
In order to analyze how the
stationarity condition changes when the set of indices changes,   the next lemma compares the generalized gradients of two approximate value functions at  points at which   the values coincide.

\begin{lem}\label{GradInclusion}

Let $U, V \subseteq Y$ be  nonempty and closed with $U\subseteq V$. If
$
\varphi_U(x)=\varphi_V(x)$ holds,
then
\begin{align*}
\partial\varphi_U(x)
\subseteq
\partial\varphi_V(x).
\end{align*}
\end{lem}

\begin{proof}
The inclusion $U \subseteq V$  implies $\varphi_U(x)\le \varphi_V(x)$, thus  the inequality
$$\varphi^\circ_U (x;d)\le \varphi^\circ_V(x;d)$$
holds
for every $d \in \mathbb{R}^n$ whenever $\varphi_V(x)=\varphi_U(x)$. This yields the  inclusion by definition of the Clarke generalized gradient.
\end{proof}

 The inclusion in Lemma~\ref{GradInclusion} implies that stationarity is preserved when  additional indices are considered, provided that the point remains feasible for  the new approximate problem.

\begin{proposition}\label{StationaryPointInclusion}

Let $U, V \subseteq Y$ be  nonempty and closed with $U\subseteq V$. If $x$ is stationary for  \eqref{VPU} and $\varphi_V(x)\le0$ holds,
then it also is stationary for  {\rm (VP$_V$)}.
\end{proposition}

\begin{proof}
Feasibility for (VP$_V$) is granted by  assumption, so we have to prove that there exists    $\lambda_V \ge 0 $ such that
\begin{align*}
&0 \in \partial f(x)+ \lambda_V \partial \varphi_V(x)+ \mathcal{N}_X(x)\\
& \lambda_V \varphi_V(x)=0.
\end{align*}
Two cases are possible. In the first we  have $\varphi_U(x)<\varphi_V(x)$ and  we also get $\varphi_U(x)<0$. Therefore,  $x$ is stationary for \eqref{VPU} along  with the multiplier $\lambda_U=0$,   i.e.,  $$0 \in  \partial f (\bar{x})+\mathcal{N}_X(\bar{x}),$$
 hence $x $ is also stationary for (VP$_V$) along with the multiplier $\lambda_V=0$.
In the case $\varphi_V(x)=\varphi_U(x)$,  we get the inclusion $$\partial \varphi_U (x) \subseteq \partial \varphi_V (x)$$ thanks to Lemma \ref{GradInclusion}, that implies the inclusion
 \begin{align*}
     &0 \in \partial f(x) + \lambda_U \partial\varphi_U(x)+ \mathcal{N}_X(x) \subseteq \partial f(\bar{x}) + \lambda_U \partial \varphi_V(\bar{x})+ \mathcal{N}_X(x)
    \end{align*}
    with $\lambda_U\varphi_V(x)=\lambda_U\varphi_U(x)=0$. Hence, $x $ is also stationary for (VP$_V$) along with the  multiplier $\lambda_V=\lambda_U$.

\end{proof}

 We now discuss the limit behavior of the approximated value function as the set of indices changes. This setting is relevant to cutting surface algorithms, in which new indices are added as
 iterations proceed.
In fact, the next result characterizes the limit behavior of the approximate value functions given by   sequence of sets $\{U_k\}_k$ in terms of the value function defined by Kuratowski superior limit $\widehat{U}$ of the sequence. Note that  the superior limit is always a nonempty compact subset of $Y$ since $Y$ is compact (see \cite[Proposition 1.1.2]{aubin_set-valued_2009}).

\begin{proposition}\label{ConvergenceGeneralized Gradients}
 Let $\{U_k\}_{k\in \mathbb{N}}$  be a   sequence of closed subsets of $Y$ with Kuratowski superior limit $\hat{U}$ and consider   $x_k \rightarrow x_\infty$.
Then, the limit superior  reads    $$\limsup_{k \rightarrow \infty}\varphi_{U_k}(x_k)=\varphi_{\hat{U}}(x_\infty).$$ Moreover,
    when the limit exists, i.e.,   $$\lim_{k \rightarrow \infty} \varphi_{U_k}(x_k)=\varphi_{\hat{U}}(x_\infty),$$ then the following implication  holds true
 $$\ \xi_k \rightarrow \xi_\infty \   with \  \xi_k \in \partial \varphi_{U_k}(x_k)  \implies \xi_\infty  \in \partial \varphi_{\hat{U}}(x_\infty).$$

 \end{proposition}

\begin{proof}
The inequality $$\limsup_{k \rightarrow \infty}\varphi_{U_k}(x_k)\le\varphi_{\hat{U}}(x_\infty)$$
follows from \cite[Lemma 1.1.6]{aubin_set-valued_2009}. Let $\hat{y}\in \arg \max_{y \in \hat{U}} g(x_\infty,y)$. The definition of Kuratowski superior limit guarantees that there exist  $y_{k_i}\in U_{k_i}$ such that $d(y_{k_i},\hat{y})\rightarrow0$. Thanks to the continuity of $g$, we obtain
$$\limsup_{k \rightarrow \infty}\varphi_{U_k}(x_k) \ge \lim_{i \rightarrow \infty }g(x_{k_i},y_{k_i})=g(x_\infty,\hat{y})=\varphi_{\hat{U}}(x_\infty)$$
by taking the limit as $i\rightarrow \infty$.

Now, assume that $ \varphi_{\hat{U}}(x_\infty) $ is the limit of the sequence $\{\varphi_{{U_k}}(x_k)\}_k$. The Clarke generalized gradient of $\varphi_{U_k}$ at $x_k$ is given by
$$\partial  \varphi_{U_k}(x_k)=\text{co}\{ \xi \in \partial_x g(x_k,y) \ : \ y \in \arg \max_{y\in U_k} g(x_k,y)\},$$
according to  Proposition \ref{ClarkeGrad}.
Since these sets are compact and convex,  each of their elements $\xi_k$ is a convex combination of at most $n+1$ points by Carathéodory's Theorem, i.e.,
$$\xi_k =\sum_{i=1}^{n+1}\gamma_i^kv_i^k,   \quad with \quad v_i^k\in \partial_xg(x_k,y^k_i) \quad and \quad \ y^k_i \in \arg \max_{y \in U_k} g(x_k,y),$$
where  $\gamma_i^k\ge 0$ and  $\gamma_1^k+\dots+\gamma_{n+1}^k =1$. Now, consider the sequences $\{\gamma^k_i\}_{k \in \mathbb{N}}$, $\{v^k_i\}_{k \in \mathbb{N}}$ and $\{y^k_i\}_{k \in \mathbb{N}}$. Since each sequence entirely lies in   a  compact set, there exists a common subsequence of indices $\{k_j\}_j$ such that
$$\gamma_i^{k_j } \rightarrow \gamma_i^\infty\in [0,1], \quad v_i^{k_j}\rightarrow v_i^\infty \in \partial_x g(x_\infty,y_i^\infty), \quad  y_i^{k_j } \rightarrow y_i^\infty \in Y$$
holds for each $i$. Notice that  the inclusion $v_i^\infty \in \partial_x g(x_\infty,y_i^\infty)$ follows from the outer semi-continuity of the set-valued  map
$\partial_x g$.
Now, we obtain both
$$\xi_\infty = \lim_{j \rightarrow \infty}\xi_{k_j}=\lim_{j \rightarrow \infty}\sum_{i=1}^{n+1}\gamma_i^{k_j}v_i^{k_j}=\sum_{i=1}^{n+1}\gamma_i^{\infty}v_i^\infty$$
and
$$\sum_{i=1}^{n+1}\gamma_i^{\infty}=\lim_{j \rightarrow \infty }\sum_{i=1}^{n+1}\gamma_i^{k_j} = 1,$$
by taking the limit as $j\rightarrow\infty$.
Moreover,  the  continuity of $g$ implies
$$g(x_\infty,y_i^\infty)=\lim_{j \rightarrow\infty}g(x_{k_j},y_i^{k_j})=\lim_{j \rightarrow \infty} \varphi_{U_{k_j}}(x_{k_j})=\lim_{k \rightarrow \infty} \varphi_{U_k}(x_k)=\varphi_{\hat{U}}(x_\infty),$$
so that  $y_i^\infty\in \arg\max_{y \in \hat{U}} g(x_\infty,y)$ and hence   $\xi_\infty\in \partial\varphi_{\hat{U}}(x_\infty)$.

\end{proof}

As a consequence, the limiting behavior of the Clarke generalized gradients can be controlled.

\begin{corollary}\label{LimitsetInclusion}

     Let $\{U_k\}_{k\in \mathbb{N}}$  be a   sequence of closed subsets of $Y$ with Kuratowski superior limit $\hat{U}$ and consider a sequence  $x_k \rightarrow x_\infty$. Then $$\Limsup_{k\rightarrow \infty}  \partial \varphi_{U_k}(x_k) \subseteq \partial \varphi_{\hat{U}}(x_\infty)$$
 \end{corollary}
\begin{proof}
 Proposition \ref{ConvergenceGeneralized Gradients} guarantees that the  limit superior  exists and is equal to $\varphi_{\hat{U}}(x_\infty)$. Hence, we can find a subsequence of indices $\{k_j\}_j$ such that   $$\lim_{j \rightarrow \infty} \varphi_{U_{k_j}}(x_{k_j})=\varphi_{\hat{U}}(x_\infty).$$
    Now, let  $\tilde{U}$ be the Kuratowski superior limit of the sequence $\{U_{k_j}\}_j$. Then,  we obtain
    $$\Limsup_{k\rightarrow \infty}  \partial \varphi_{U_k}(x_k) \subseteq \partial \varphi_{\tilde{U}}(x_\infty) \subseteq \partial \varphi_{\hat{U}}(x_\infty),$$
    where the first inclusion holds
  thanks to  Proposition  \ref{ConvergenceGeneralized Gradients}, and the second  inclusion
    holds
    thanks to  Lemma \ref{GradInclusion} since $\tilde{U} \subseteq\hat{U}$ and $\varphi_{\tilde{U}}(x_\infty)=\varphi_{\hat{U}}(x_\infty)$.
\end{proof}

The above results describe the convergence of the approximate value functions and their generalized gradients as the set of indices changes.
We now study how stationary points for (VP$_{U_k}$) behave when $U_k$ changes,  showing that  stationarity is preserved in the limit whenever  the  associated sequence of multipliers  remains bounded.

\begin{proposition}\label{StationaryPointLimit}
      Let $\{U_k\}_{k\in \mathbb{N}}$  be a   sequence of closed subsets of $Y$ with Kuratowski superior limit $\hat{U}$ and consider   $x_k \rightarrow x_\infty$ such that $x_k$ is stationary for (VP$_{U_k}$) along with the multiplier $\lambda_k$. If the sequence  $\{\lambda_k\}_k$ is bounded,  then $x_\infty$ is stationary for (VP$_{\hat{U}}$).
\end{proposition}

\begin{proof}
  Proposition \ref{ConvergenceGeneralized Gradients} guarantees
$$\limsup_{k\rightarrow\infty}\varphi_{U_k}(x_k)=\varphi_{\hat{U}(x_\infty)},$$
    that  implies  $\varphi_{\hat{U}}(x_\infty)\le 0$ since   $\varphi_{U_k}(x_k)\le0$ holds for every $k$. By taking the Kuratowski superior limit of the stationarity conditions for (VP$_{U_k}$), we get:
    \begin{align*}
        0 &\in \Limsup_{k\rightarrow \infty}(\partial f (x_k)+\lambda_k \partial \varphi_{U_k}(x_k)+ \mathcal{N}_X(x_k))\\
        & \subseteq \Limsup_{k\rightarrow \infty}(\partial f (x_k))+\Limsup_{k\rightarrow \infty}(\{\lambda_k\}) \Limsup_{k\rightarrow \infty}(\partial \varphi_{U_k}(x_k))+ \Limsup_{k\rightarrow \infty}(\mathcal{N}_X(x_k))\\
        &\subseteq
        \partial f (x_\infty)+\Lambda_\infty \partial \varphi_{\hat{U}}(x_\infty)+ \mathcal{N}_X(x_\infty)
    \end{align*}
    where $\Lambda_\infty$ is the Kuratowski superior limit of the sequence of the multipliers $\{\lambda_k\}_k$.
    Since the sequences of indices $\{\lambda_k\}_k$  and $\{\varphi_{U_k}(x_k)\}_k$  are bounded, there exists a subsequence $\{k_j\}_j$ such that both the limits
    $$\lambda_{k_j}\rightarrow \lambda_\infty  \quad and \quad \lim_{k_j\rightarrow\infty}\varphi_{U_{k_j}}(x_{k_j})=\varphi_{\hat{U}(x_\infty)}$$
    hold, where $\lambda_\infty \in \Lambda_\infty$. Taking the limit as $j \rightarrow \infty$ in the complementarity conditions $$\lambda_{k_j}\varphi_{U_{k_j}}(x_{k_j})=0$$
    yields that $\lambda_\infty\varphi_{\hat{U}}(x_\infty)=0$, hence $x_\infty$ is stationary for  (VP$_{\hat{U}}$).
\end{proof}

\subsection{Stationarity under penalization}

Instead of computing stationary points of \eqref{VPU} directly, we incorporate the violation of the constraint into the objective function. Consider the positive part of the approximate value function, i.e.,
\begin{align*}
\varphi_U^+(x)
&:=
\max\{0,\varphi_U(x)\}
\end{align*}
in order to measure the
violation of the constraint. Indeed, it vanishes at any feasible point of \eqref{VPU} while  being positive elsewhere.

Given a penalty parameter $\rho>0$,  consider the penalized problem
\begin{align}
\min_{x\in X}\quad
& f(x)+\rho\varphi_U^+(x).
\tag{VP$_U^\rho$}\label{PenalizedVPU}
\end{align}

Since $\varphi_U$ is Clarke regular, Theorem 2.3.9 in \cite{clarke_optimization_1990} guarantees that  $\varphi_U^+$ is also Clarke regular and its generalized gradient at any  point $x \in X$ is given by
 \begin{align*}
    \partial  \varphi_U^+ (x):=
    \left\{
        \begin{array}{ll}
           \{0\}  &  {\text{if}} \quad \varphi_U(x)<0\\
             \{\theta \xi \ | \  \theta \in [0,1],\xi\in \partial \varphi_U(x)  \} , \quad   &  {\text{if}} \quad \varphi_U(x)=0\\
            \partial \varphi_U(x) &  {\text{if}} \quad \varphi_U(x)>0.
        \end{array}
        \right.
    \end{align*}

 Since the constraint violation is incorporated into the objective function, the stationarity condition for  \eqref{PenalizedVPU} does not include a complementarity condition. A point $x\in X$ is stationary for  \eqref{PenalizedVPU} if it satisfies the inclusion:
\begin{align*}
0
\in
\partial f(x)
+
\rho\partial\varphi_U^+(x)
+
\mathcal{N}_X(x).
\end{align*}
 Differently from \eqref{VPU}, stationarity is satisfied at any   local minimizer without  the need of any constraint qualification. Moreover, when a stationary point for  \eqref{PenalizedVPU} is feasible for \eqref{VPU}, then it is also stationary for the latter, and  the  penalty parameter  yields an upper bound on the magnitude of  at least one  multiplier.

\begin{proposition}\label{PenaltyFeasibility}

Let $x\in X$ be a stationary point of problem \eqref{PenalizedVPU}. If $x$ is feasible for  \eqref{VPU}, then there exists $\lambda_U\le \rho$ such that $x$ is stationary for  $\eqref{VPU}$ along with multiplier $\lambda_U$.
\end{proposition}

\begin{proof}
Since $x$ is feasible for \eqref{VPU},  two cases  may occur. If $\varphi_U(x)<0$ holds, then $\partial  \varphi_U^+ (x)={0}$, so that the point is stationary for \eqref{VPU} with multiplier equal to $0$. Otherwise, if $\varphi_U(x)=0$, then there exist $\theta \in [0,1]$ and $\xi \in \partial \varphi_U(x)$ such that
\begin{align*}
&0 \in \partial f(x) + \rho \theta\xi+ \mathcal{N}_X(x)
\end{align*}
 since $x$ is stationary for \eqref{PenalizedVPU}. Hence, $x$ is stationary for  \eqref{VPU} along with the  multiplier $\lambda_U=\rho\theta\le \rho$.
\end{proof}

\section{Algorithms}

\label{sec:algorithms}

 Algorithms that exploit the approximate problems  introduced in the previous section are devised.  All the algorithms  follow the same general scheme: at each iteration, a stationary point of a penalized approximation is computed, while the penalty parameter is adjusted to promote feasibility for the corresponding constrained approximation. Cutting surfaces are used to progressively improve the approximation of  the constraints of the original semi-infinite problem \eqref{sip}. In this way, the penalty mechanism controls feasibility with respect to the current approximation, while the cutting surface updates drive the iterates toward feasibility for the original problem. A
 basic penalty algorithm   is devised according to this strategy and  proximal methods are designed by also adding a quadratic regularization term to the objective function. The first proximal algorithm directly exploits $f$, whereas the second replaces it with its first-order linearization.  A further variant allows discarding the  indices that correspond to the constraints that are inactive at the current iterate. All the methods converge to stationary points for \eqref{sip} when the sequence of penalty parameters remains bounded. Finally, an exact penalty theorem is established provided  that eNNAMCQ holds at each feasible and infeasible point in $X$.

\subsection{Basic penalty algorithm}
\label{subsec:prototype}

We first introduce the
basic penalty algorithm that directly combines the cutting surface and penalty techniques. Starting from a nonempty  index set $Y_0$ and an initial penalty parameter $\rho_0$, each iteration computes a stationary point for  (VP$_{Y_k}^{\rho_k}$). If this point is feasible  for ~\eqref{sip}, the algorithm terminates since the point is also stationary for \eqref{sip}. Indeed, Proposition \eqref{PenaltyFeasibility} guarantees that the point is stationary for (VP$_{Y_k}$) and Proposition \eqref{StationaryPointInclusion} guarantees that it is also stationary for \eqref{sip}. If the point is infeasible also  for (VP$_{Y_k}$), the penalty parameter is increased to promote feasibility at the next iterations. Then, the approximation is refined by adding  a new index to $Y_k$ in order to exclude the current iterate from successive iterations.

\begin{algorithm}[h]
\caption{Basic penalty}\label{alg1}
\begin{algorithmic}[1]
\State \textbf{input:}  $Y_0 \subset Y$, $\rho_0>0$, $\sigma>1$.
\For{$k=0,1,\dots$}
\State Compute $x_k$ stationary for problem:
$ \min\{f(x) + \rho_k \varphi^+_{ \ Y_k}(x)\ : \ x \in X \ \}$
\If{$\varphi(x_k)\le0$}
\State {\textbf{stop}}

\Else
\If{$\varphi^+_{Y_k} (x_k)>0$}
    \State{$ \rho_{k+1}=\sigma \rho_k$}
\Else
\State $\rho_{k+1}=\rho_k$
\EndIf
\State Compute $y_{k+1}$ such that:
			$y_{k+1} \in \arg\max \{ \ g(x_k,y) \ : \ y \in Y \ \}$
\State $Y_{k+1}=Y_k \cup \{y_{k+1}\}$

\EndIf
\EndFor
\end{algorithmic}
\end{algorithm}

The convergence of the algorithm depends on the behavior of the sequence of the  penalty parameters. Indeed, the method either terminates after finitely many iterations or generates an infinite sequence whose limit $\rho_\infty$ determines the kind of convergence. Indeed,
 every cluster point of the iterate sequence is stationary for \eqref{sip} when  this limit is finite, whereas  stationarity is guaranteed only for cluster points satisfying eNNAMCQ when the limit is infinite.

\begin{theorem}\label{ThmAlg1}
    If Algorithm \ref{alg1} terminates in a finite number of   iterations $\bar{k}$, then $x_{\bar{k}}$ is stationary for \eqref{sip}. Otherwise,

    the following hold true:
    \begin{enumerate}
        \item if $\rho_\infty<+\infty$, then every cluster point of $\{x_k\}_{k\in \mathbb{N}}$ is a stationary point for  \eqref{sip}.
        \item if $\rho_\infty=+\infty$, then every cluster point  of $\{x_k\}_{k\in \mathbb{N}}$ that satisfies eNNAMCQ is also  a stationary point for  \eqref{sip}.
\end{enumerate}
\end{theorem}
\begin{proof}

If Algorithm \ref{alg1} stops at iteration $\bar{k}$, then  $x_{\bar{k}}$ is stationary for problem (VP$_{Y_k}^{\rho_k}$).

 The inequality $\varphi_{Y_k}(x_{\bar{k}}) \le 0$ holds thanks to  the stopping criterion, hence $x_{\bar{k}}$ is also stationary for (VP$_{Y_k}$)  by Proposition \ref{PenaltyFeasibility}. Moreover, Proposition \ref{StationaryPointInclusion} implies that  $x_{\bar{k}}$ is stationary for (SIP)  since $\varphi(x_{\bar{k}}) \le 0$.

Now,  assume that Algorithm \ref{alg1} never stops and  let $x_\infty$ be a cluster point of the sequence $\{x_k\}_k$. Then, there exists a subsequence of indices $\{k_j\}_j$ such that $x_{k_j}\rightarrow x_\infty$. Thus, the whole sequence $\{x_{k_j}\}_j$ lies in a compact set, so that  $(A3)$ implies that there exists   $H>0$  such that the inequality
$$||g(x_{{k_j}},y)- g(x_{k_{j-1}},y)||\le H||x_{{k_j}}-x_{k_{j-1}}||$$
holds for every $y \in Y$. Hence, the following chain of inequalities
holds
\begin{align*}
    -H||x_{{k_j}}-x_{k_{j-1}}||&\le g(x_{{k_j}},y_{k_{j-1}})-g(x_{k_{j-1}},y_{k_{j-1}})\\
    &\le \varphi_{Y_{k_j}}(x_{{k_j}})- \varphi_{Y_{k_j}}(x_{k_{j-1}})\\
    &\le g(x_{{k_{j}}},y_{k_{j}})-g(x_{k_{j-1}},y_{k_{j}})\\
    &\le H||x_{{k_j}}-x_{k_{j-1}}||,
\end{align*}
where the first and last inequalities follow from the Lipschitz bound, while the second and the third follows  both from the equalities  $$ \varphi_{Y_{k_j}}(x_{{k_j}})=g(x_{k_j},y_{k_{j}}) \ and \ \varphi_{Y_{k_{j}}}(x_{{k_{j-1}}})=g(x_{k_{j-1}},y_{k_{j-1}})$$ and the inequalities $$g(x_{k_j},y_{k_{j-1}})\le\varphi_{Y_{k_j}}(x_{{k_j}})\  and  \ g(x_{k_{j-1}},y_{k_{j}})\le\varphi_{Y_{k_j}}(x_{{k_{j-1}}}),$$
which derive from  steps 12 and 13 of the algorithm and the definition of the approximate value function.
Hence, we obtain
$$|\varphi_{Y_{k_j}}(x_{{k_j}})- \varphi_{Y_{k_j}}(x_{k_{j-1}})|\le H||x_{{k_j}}-x_{k_{j-1}}||,$$
which yields
the following inequalities $$\varphi(x_{k_j}) \ge \varphi_{Y_{k_j}}(x_{{k_j}})\ge \varphi_{Y_{k_j}}(x_{k_{j-1}})-H||x_{{k_j}}-x_{k_{j-1}}||=\varphi(x_{k_{j-1}})-H||x_{k_j}-x_{k_{j-1}}||$$
since $\varphi\ge \varphi_{Y_{k_j}}$ holds by   the definition of the value function and $\varphi(x_{j-1})= \varphi_{Y_{k_j}}(x_{j-1})$ holds by step 12 of the algorithm.
Now, let $\hat{Y}$ be  the Kuratowski superior limit of the sequence $\{Y_{k_j}\}_j$.
Taking the limit as $j \rightarrow \infty$ yields $$\lim _{j \rightarrow \infty}\varphi_{Y_{k_j}}(x_{k_j})=\varphi(x_\infty),$$
that implies $\varphi_{\hat{Y}}(x_\infty)=\varphi(x_\infty)$ thanks to Proposition \ref{ConvergenceGeneralized Gradients}. Hence, we obtain the inclusion $$\partial \varphi_{\hat{Y}}(x_\infty)
    \subseteq \partial \varphi(x_\infty) .$$

1. Since   the sequence  $\{\rho_{k_j}\}_j$ is bounded, then there exists a $\bar{j}$ such that  $\varphi_{Y_{k_j}} (x_{k_j})\le0$  for every $j \ge \bar{j}$. Thus,  every $x_{k_j}$ is also stationary for (VP$_{Y_{k_j}}$) along with  some  multiplier  $\lambda_{k_j}\le \rho_{k_j}$ when $j$ is sufficiently large. The sequence of the multipliers is bounded by $\rho_\infty$, so that  $x_\infty$ is stationary for (VP$_{\hat{Y}}$) by  Proposition \ref{LimitsetInclusion}. Moreover, we have
$$0\ge \lim _{j \rightarrow \infty}\varphi_{Y_{k_j}}(x_{k_j})=\varphi(x_\infty),$$
 so that  $x_\infty$ is feasible for \eqref{sip} and therefore stationary    by Proposition \ref{StationaryPointInclusion}.

2. Now, assume that  the sequence  $\{\rho_{k_j}\}_j$ is not bounded and that the cluster point $x_\infty$ satisfies  eNNAMCQ. Then, two possible cases arise:
\begin{itemize}
    \item[a)] There exists $\bar{j}$ such that   $\varphi_{Y_{k_j}} (x_{k_j})\le0$  for every $j \ge \bar{j}$,
    \item[b)] $\varphi_{Y_{k_j}}(x_{k_j})>0$ happens infinitely often.
\end{itemize}

When a) holds true, then every $x_{k_j}$ is also stationary for (VP$_{Y_{k_j}}$) along with some  multiplier  $\lambda_{k_j}\le \rho_{k_j}$ when $j$ is sufficiently large. If the sequence of the multipliers is bounded, then we can argue as in previous situation (point 1.) to obtain that $x_\infty $ is  stationary for \eqref{sip}. Ab absurdo, assume that sequence of the multipliers is not bounded. Since  $x_{k_j}$ is stationary for (VP$_{Y_{k_j}}$) we get   $$0 \in  \dfrac{\partial f(x_{k_j})}{\lambda_{k_j}} +  \partial \varphi_{Y_{k_j}}(x_{k_j}) + \mathcal{N}_X(x_{k_j}).$$
 By considering the Kuratowski superior limit we obtain:
\begin{align*}
    0 & \in \Limsup_{j\rightarrow \infty}({\partial f(x_{k_j})}/{\lambda_{k_j}} +  \partial \varphi_{Y_{k_j}}(x_{k_j}) + \mathcal{N}_X(x_{k_j})) \\
   & \subseteq \Limsup_{j\rightarrow \infty}({\partial f(x_{k_j})}/{\lambda_{k_j}}) +  \Limsup_{j\rightarrow \infty}(\partial \varphi_{Y_{k_j}}(x_{k_j})) + \Limsup_{j\rightarrow \infty}(\mathcal{N}_X(x_{k_j})) \\
   &\subseteq \partial \varphi_{\hat{Y}}(x_\infty) + \mathcal{N}_X(x_\infty)
    \subseteq \partial \varphi(x_\infty) + \mathcal{N}_X(x_\infty).
\end{align*}
Moreover,  we have $\varphi(x_{k_j})>0$ since the stopping criterion is never satisfied, which implies that   $\varphi(x_\infty) \ge 0$, in contradiction with  eNNAMCQ.

Assume b) holds true. Note that $$0 \in  \dfrac{\partial f(x_{k_j})}{\rho_{k_j}} +  \partial \varphi_{Y_{k_j}}(x_{k_j}) + \mathcal{N}_X(x_{k_j})$$
whenever $\varphi_{k_j}(x_{k_j}) > 0$ since each $x_{k_j}$ is stationary for (VP$_{Y_{k_j}}^{\rho_{k_j}}$).  As a consequence, the
  Kuratowski superior limit of the above sequence of sets satisfies
\begin{align*}
    0 & \in \Limsup_{j\rightarrow \infty}({\partial f(x_{k_j})}/{\rho_{k_j}} +  \partial \varphi_{Y_{k_j}}(x_{k_j}) + \mathcal{N}_X(x_{k_j})) \\
   & \subseteq \Limsup_{j\rightarrow \infty}({\partial f(x_{k_j})}/{\rho_{k_j}}) +  \Limsup_{j\rightarrow \infty}(\partial \varphi_{Y_{k_j}}(x_{k_j})) + \Limsup_{j\rightarrow \infty}(\mathcal{N}_X(x_{k_j})) \\
   &\subseteq \partial \varphi_{\hat{Y}}(x_\infty) + \mathcal{N}_X(x_\infty)
    \subseteq \partial \varphi(x_\infty) + \mathcal{N}_X(x_\infty).
\end{align*}
Moreover,  we have $\varphi(x_{k_j})>0$ since the stopping criterion is never satisfied, that implies  $\varphi(x_\infty) \ge 0$ in contradiction with  eNNAMCQ.
\end{proof}

\subsection{Proximal methods}
\label{subsec:proximal-methods}

Algorithm~\ref{alg1} requires the computation of a stationary point of a penalized approximation of the semi-infinite problem at each iteration. As new constraints are added, the feasible region changes and the stationary point may move significantly between iterations, causing  the sequence of iterates to jump around. To reduce this possible oscillatory behavior, we introduce a proximal variant of the basic penalty algorithm by adding a quadratic regularization term centered at the current iterate.
At each iteration $k$, an approximate solution of the problem
\begin{align}\label{regularizedprob}
\min_{x\in X} \quad
& f(x)+\rho_k\varphi^+(x)
+\frac{1}{2\gamma_k}||x-x_k||^2
\end{align}
is computed, where the proximal parameter $\gamma_k$ controls the relevance of the regularization.
To this aim, an inner loop indexed by $i$ is performed. At each inner iteration, the next iterate $z_{i,k}$ is computed by solving  an approximation of   \eqref{regularizedprob} where  $\varphi^+$ is replaced with  $\varphi^+_{Y_{i,k}}$.  The approximation is then refined by adding   a new index $y_{i,k}$ to  $Y_{i,k}$, that corresponds to the constraint that the current iterate violates the most. The inner loop terminates when the approximate value function is sufficiently close to the exact one at $z_{i,k}$ according to the prescribed tolerance $\delta_k$, i.e.,
$$\varphi^+(z_{i,k})-\varphi_{Y_{i,k}}^+(z_{i,k}) \le \delta_k.$$
When this happens,  the final inner iterate  becomes the new center of the proximal term, while the two indices  $w_{k+1}$ and  $y_{k+1}$ that satisfy
$$\varphi(z_{i,k})=g(z_{i,k},y_{k+1}) \quad and \quad \varphi_{Y_{i,k}}(z_{i,k})=g(z_{i,k},w_{k+1}) $$
 are added to the index set $Y_{k}$. While the index  $y_{k+1}$ allows excluding the current iterate  from the next iterations, $w_{k+1}$ preserves the latest information on  the approximate value function $\varphi_{Y_{i,k}}$  without carrying over the entire index set $Y_{i,k}$.
Finally, the penalty parameter is increased  whenever the iterate is infeasible for the approximate problem (VP$_{Y_{i,k}}$) to promote feasibility at the next iterations.

\begin{algorithm}[h]
\caption{Penalty with proximal regularization}
\label{alg3}
\begin{algorithmic}[1]
\State \textbf{Input:}
$x_0\in X$,  $Y_0\subseteq Y$, $\rho_0>0$, $\sigma>1$ and
$\{\delta_k\}_k,\{\gamma_k\}_k \subseteq \mathbb{R}_+$
\For{$k=0,1,\ldots$}
\State Set $Y_{0,k}:=Y_k$
\For{$i=0,1,\ldots$}
\State Compute
$
z_{i,k}\in
\arg\min
\{
f(x)
+
\rho_k\varphi_{Y_{i,k}}^+(x)
+
\dfrac{1}{2\gamma_k}||x-x_k||^2
\ : \ x \in X\}
$
\State Compute
$
y_{i,k}\in
\arg\max\{g(z_{i,k},y) \ : \ y \in Y\}
$
\State Compute
$
w_{i,k}\in
\arg\max\{g(z_{i,k},y) \ : \ y \in Y_{i,k}\}
$
\If{
$\varphi^+(z_{i,k})-\varphi_{Y_{i,k}}^+(z_{i,k}) \le \delta_k$

}
\State
$x_{k+1}:= z_{i,k}$,
$y_{k+1}:= y_{i,k}$,
and
$w_{k+1}:= w_{i,k}$
\State
$Y_{k+1}:=
Y_k\cup\{w_{k+1},y_{k+1}\}$
\State \textbf{break}
\Else
\State
$Y_{i+1,k}:=
Y_{i,k}\cup\{y_{i,k}\}$
\EndIf
\EndFor
\If{$g(x_{k+1},w_{k+1})>0$}
\State $\rho_{k+1}:=\sigma\rho_k$
\Else
\State $\rho_{k+1}:=\rho_k$
\EndIf
\If{
$x_{k+1}=x_k$
and
$\varphi(x_{k+1})\leq 0$
}
\State \textbf{stop}
\EndIf
\EndFor
\end{algorithmic}
\end{algorithm}

Beyond the feasibility of the iterate, i.e., $\varphi(x_{k+1})\le0$,
the stopping criterion  additionally requires $x_{k+1}=x_k$ due to the presence of the quadratic term in the objective function. Indeed,  the first-order optimality condition of problem at line 5 read
$$ \dfrac{x_{k}-x_{k+1}}{\gamma_k}
\in
\partial f(x_{k+1})
+
\rho_k\partial\varphi^+_{Y_{i,k}}(x_{k+1})
+
\mathcal{N}_X(x_{k+1}),$$ so that $x_{k+1}$ is stationary for (VP$^{\rho_k}_{Y_{i,k}} $) whenever $x_{k+1}=x_{k}$ holds.
Therefore, $x_{k+1}$ is a stationary point for \eqref{sip} thanks to Proposition \ref{StationaryPointInclusion} and \ref{PenaltyFeasibility}.

The analysis requires the following quadratic lower bound on the objective function, which is  always satisfied by weak convex functions.
\begin{itemize}
    \item[\textbf{(A4)}] There exist $M\in\mathbb{R}$ and $\mu\geq 0$ such that
\begin{align*}
f(x)
\geq
M-\frac{\mu}{2}||x||^2
\qquad
\forall x\in X.
\end{align*}
\end{itemize}
This assumption ensures that the  objective function of problem at line $5$ of the algorithm  is coercive when the proximal parameter is sufficiently small, so that the inner loop terminates in a finite number of steps.

A  variant of Algorithm~\ref{alg3} can be employed when $f$ is $L$-smooth, that is the following condition holds:

\begin{itemize}
    \item[\textbf{(A5)}]  $f$ is continuously differentiable  and its gradient is Lipschitz continuous on $X$, i.e., there exists  $L>0$ such that
\begin{align*}
||\nabla f(x)-\nabla f(z)||
\leq
L||x-z||
\end{align*}
holds for every $x,z\in X$.
\end{itemize}

The above assumption  implies  that the lower bound
$$f(x)\ge
f(z)+\langle \nabla f(z),x-z\rangle
-\frac{L}{2}||x-z||^2,$$
holds for every $x,z \in X$ (see, \cite[Lemma 1.2.3]{Nesterov2004}). Therefore,   $f$ also satisfies  (A4)  since, for any fixed $\bar{z} \in X$, the right hand side of the inequality
$$f(x)+\frac{\mu}{2}||x||^2\ge
f(\bar{z})+\langle \nabla f(\bar{z}),x-\bar{z}\rangle
-\frac{L}{2}||x-\bar{z}||^2+\frac{\mu}{2}||x||^2$$
is a strongly convex function in $x$ whenever $\mu>L$. Moreover,
when (A5) holds  true,   the objective function in the proximal subproblem can be  replaced with  its  linearized first-order approximation.

\medskip

\textbf{Algorithm 3:}
step $5$ in Algorithm \ref{alg3} is replaced with
$$ 5': \text{Compute } z_{i,k} \in
 \arg \min\{\nabla f(x_k)^Tx+\rho_k \varphi^+_{Y_{i,k}}(x) + \dfrac{1}{2\gamma_k}{||x-x_k
 ||}^2\ : \ x \in X\}$$

To ensure convergence, the  parameters of the algorithms are chosen according to the following criteria:
\begin{itemize}

\item the sequence of tolerances $\{\delta_k\}_k$ is summable,

\item the  sequence of proximal parameters $\{\gamma_k\}_k$ satisfies 
\begin{align*}
    \gamma_k\in
    \left[
        \bar{\gamma},
        \bar{\beta}
    \right]
\end{align*}
{for some  given $\bar{\beta},\bar{\gamma}$ such that $0<\bar{\gamma} < \bar{\beta}< 1/\mu,$
in Algorithm~2,} and
\begin{align*}
    \gamma_k\in
    \left[
        \bar{\gamma},
        \bar{\beta}
    \right]
\end{align*}
for some  given $\bar{\beta},\bar{\gamma}$ such that $0<\bar{\gamma} < \bar{\beta}< 1/L,$ in Algorithm~3,
\end{itemize}
together with  an initial index set $Y_0$ such that $f$ is bounded from below on
\begin{align*}
\left\{
x\in X
\mid
g(x,y)\leq 0,\ \forall y\in Y_0
\right\}.
\end{align*}

The next result shows that the inner loop  terminates in a finite number of iterations for any given level of accuracy $\delta_k>0$.

\begin{theorem}
    Assume  (A4) for Algorithm 2 while (A5) for Algorithm 3.  Then, the inner loop between lines 4-15 terminates in a finite number of iterations.
\end{theorem}

\begin{proof}
    Consider Algorithm 2 and let the iteration $k$ be fixed. Any inner loop $i$ satisfies
    \begin{align*}
   M - \dfrac{\mu}{2}||z_{i,k}||^2 + \dfrac{1}{2\gamma_k}{||z_{i,k}-x_k||}^2&\le  f(z_{i,k})  + \dfrac{1}{2\gamma_k}{||z_{i,k}-x_k||}^2\\
   &\le
        f(z_{i,k}) + \rho_k \varphi^+_{Y_{i,k}}(z_{i,k}) + \dfrac{1}{2\gamma_k}{||z_{i,k}-x_k||}^2 \\
        &\le  f(x_k)+ \rho_k \varphi^+_{Y_{i,k}}(x_k)\\
        & \le f(x_k)+ \rho_k \varphi^+(x_k),
    \end{align*}
    where the first inequality follows from (A4), the second from the non-negativity of $\varphi^+_{Y_{i,k}}$, the third from the fact that $z_{i,k}$ is a minimizer of the problem at line 5 and the fourth one from $\varphi^+_{Y_{i,k}}\le \varphi^+$.
    Since    $1/\gamma_k>\mu$,  the function $$h(z)= M - \dfrac{\mu}{2}||z||^2 + \dfrac{1}{2\gamma_k}{||z-x_k||}^2$$
    is coercive. Since the  quantity $f(x_k)+ \rho_k \varphi^+(x_k)$ is an upper bound for  the  sequence $\{h(z_{i,k})\}_i$,  the sequence $\{z_{i,k}\}_{i}$ is bounded. Hence,  there exists a subsequence of indices $\{i_j\}_j$ such that  $z_{{i_j},k}\rightarrow z_{\infty,k}$ for some suitable $z_{\infty,k} \in X$. Exploiting the same arguments of the  proof of Theorem \ref{ThmAlg1}, there exits a $H>0$ such that $$||\varphi_{Y_{i_j,k}}(z_{{i_j,}})- \varphi_{Y_{i_j,k}}(z_{i_{j-1},k})||\le H||z_{{i_j,k}}-z_{i_{j-1},k}||$$
    holds for every $j$. Moreover, since $\varphi_{Y_{i_j,k}}(z_{i_{j-1},k})= \varphi(z_{i_{j-1},k})$ holds thanks to steps 6 and 13 of the algorithm, we obtain
    \begin{align*}
        \varphi(z_{i_j,k})\ge \varphi_{Y_{i_j,k}}(z_{i_j,k})&\ge \varphi_{Y_{i_j,k}}(z_{i_{j-1},k})-H||z_{i_{j-1},k}-z_{i_j,k}||\\
        &=\varphi(z_{i_{j-1},k})-H||z_{i_{j-1},k}-z_{i_j,k}||.
    \end{align*}
   Therefore, $\varphi_{Y_{i_j,k}}(z_{i_j,k}) \rightarrow \varphi(z_{\infty,k})$ follows by taking the limit as $j \rightarrow \infty$. Thus,
    $$\varphi^+(z_{i_j,k})-\varphi_{Y_{i_j,k}}^+(z_{i_j,k}) \le \delta_k$$
    is  satisfied whenever $i_j$ is  large enough.

    Now, consider   Algorithm 3 and let iteration $k$ be fixed. Any inner loop $i$ satisfies  \begin{align*}
    \nabla f(x_k)^Tz_{i,k}  + \dfrac{1}{2\gamma_k}{||z_{i,k}-x_k||}^2
   &\le
        f(x_k)^Tz_{i,k} + \rho_k \varphi^+_{Y_{i,k}}(z_{i,k}) + \dfrac{1}{2\gamma_k}{||z_{i,k}-x_k||}^2 \\
        &\le  \nabla f(x_k)+ \rho_k \varphi^+_{Y_{i,k}}(x_k) \le \nabla f(x_k)+ \rho_k \varphi^+(x_k),
    \end{align*}
   where the first inequality follows from the non-negativity of  $\varphi^+_{Y_{i,k}}$, the second from the fact that $z_{i,k}$ is a minimizer of the problem at line $5'$ and the third  from $\varphi^+_{Y_{i,k}}\le \varphi^+$.
    Since the function $$h(z)=\nabla f(x_k)^Tz+ \dfrac{1}{2\gamma_k}{||z-x_k||}^2$$
    is coercive, we can exploit the same arguments of Algorithm 2 to conclude the proof.
    \end{proof}

\begin{theorem}
    Assume  (A4) for Algorithm 2 while (A5) for Algorithm 3.   If the sequence of penalty parameters
     is bounded, then every cluster point of the sequence of iterates is stationary  for \eqref{sip}.
\end{theorem}

 \begin{proof}
    We only show the convergence for Algorithm 3, as the proof can be easily adjusted for Algorithm 2. First, we prove that $||x_{k+1}-x_k||\rightarrow 0$. To this aim,  consider the  set $$W_k:=Y_k\cup \{w_{k+1}\}$$
    and the set $Y_{i_k,k}$ obtained  at the last iteration $i_k$ of the inner loop of the outer iteration $k$.
    Notice that  $x_{k+1}=z_{i_k,k}$ holds 
    thanks to step 9 of the algorithm, hence $x_{k+1}$ is a minimizer for the problem
     \begin{align}
         \min_{x \in X} \nabla f(x_{k})^Tx+ \rho_k \varphi^+_{Y_{i_k,k}}(x)+\dfrac{1}{2\gamma_k}||x-x_k||^2.\label{proxlinprob}
     \end{align}
     The following chain  holds
     \begin{align}
     \rho_k \varphi^+_{W_k}(x_{k+1})+&\nabla f(x_{k})^Tx_{k+1}+\dfrac{1}{2\gamma_k}||x_{k+1}-x_k||^2
=\notag\\
&=\rho_k \varphi^+_{Y_{i_k,k}}(x_{k+1})+\nabla f(x_{k})^Tx_{k+1}+\dfrac{1}{2\gamma_k}||x_{k+1}-x_k||^2\notag\\
     & \le \rho_k \varphi^+_{Y_{i_k,k}}(x_{k})+\nabla f(x_{k})^Tx_{k}\notag\\&\le\rho_k \varphi^+(x_k)+\nabla f(x_{k})^Tx_{k},\notag
     \end{align}
     where the  equality follows from $g(x_{k+1},w_{k+1})=\varphi_{Y_{i_k,k}}(x_{k+1})$, the first  inequality  from the fact that $x_{k+1}$ is a minimizer of problem \eqref{proxlinprob} and the second inequality  from $\varphi^+_{Y_{i,k}}\le \varphi^+$.
     Thus, we obtain
     \begin{align}
         \nabla f(x_{k})^T(x_{k+1}-x_{k})+\dfrac{1}{2\gamma_k}||x_{k+1}-x_k||^2\le\rho_k (\varphi^+(x_k)-\rho_k \varphi^+_{W_k}(x_{k+1})) \label{ineq1}
     \end{align}
     by rearranging the terms.
      Assumption (A5) guarantees that
$$\nabla f(x_{k})^T(x_k-x_{k+1})-\dfrac{L}{2}||x_{k+1}-x_k||^2\le f(x_k)-f(x_{k+1}) $$
holds by \cite[Lemma 1.2.3]{Nesterov2004}.
By summing  the above inequality with \eqref{ineq1}, we obtain
    \begin{align*}
    \dfrac{\sigma_k}{2}||x_{k+1}-x_k||^2\le f(x_{k})-f(x_{k+1}) + \rho_k(\varphi^+(x_k)-\varphi_{W_k}^+(x_{k+1}))
\end{align*}
where $\sigma_k=1/\gamma_k-L$, which leads to the inequality
\begin{align}\label{sumineq}
    \sum_{k=0}^N\dfrac{\sigma_k}{2}||x_{k+1}-x_k||^2\le f(x_0)-f(x_{N+1}) + \sum_{k=0}^N\rho_k(\varphi^+(x_k)-\varphi_{W_k}^+(x_{k+1})).
\end{align}
     Since the sequence $\{\rho_k\}_k$ is bounded, there exists $\bar{k}$ such that $\varphi_{W_k}(x_{k+1})\le 0$ and  $\rho_k=\rho_{\bar{k}}$ hold for every $k \ge \bar{k}$. Moreover, for any $N\ge \bar{k}$,  we obtain
     \begin{align*}
         &\sum_{k=0}^N\rho_k(\varphi^+(x_k)-\varphi_{W_k}^+(x_{k+1}))=\\
         &=\sum_{k=0}^{\bar{k}-1}\rho_k(\varphi^+(x_k)-\varphi_{W_k}^+(x_{k+1}))+\rho_{\bar{k}}\sum_{k=\bar{k}}^{N}(\varphi^+(x_k)-\varphi_{W_k}^+(x_{k+1}))\\
         &=S+\rho_{\bar{k}}\sum_{k=\bar{k}}^{N}(\varphi^+(x_k)-\varphi^+(x_{k+1})+\varphi^+(x_{k+1})-\varphi_{W_k}^+(x_{k+1}))\\
         &=S+\rho_{\bar{k}}(\varphi^+(x_{\bar{k}})-\varphi^+(x_{N+1}))+\rho_{\bar{k}}\sum_{k=\bar{k}}^{N}(\varphi^+(x_{k+1})-\varphi_{W_k}^+(x_{k+1}))\\
         &\le S+\rho_{\bar{k}}\varphi^+(x_{\bar{k}})+\rho_{\bar{k}}\sum_{k=\bar{k}}^{N}\delta_k,
     \end{align*}
     where $$S=\sum_{k=0}^{\bar{k}-1}\rho_k(\varphi^+(x_k)-\varphi_{W_k}^+(x_{k+1}))$$ is a  quantity that  only depends on $\bar{k}$. Also,  we have
     \begin{align*}
     f(x_{N+1})&\ge \min \{ \ f(x) \ | \ x \in X,\  \varphi_{W_N}(x)\le 0 \  \}\\
     &\ge\min \{ \ f(x) \ | \ x \in X,\  \varphi_{W_0}(x)\le 0 \  \}\\
     &>-\infty
     \end{align*}
     since $\varphi_{W_N}(x_{N+1})\le 0$.
     Hence, the sequence  $\{f(x_{N+1})\}_{N\ge\bar{k}}$ is lower bounded.
     Taking the limit as $N \rightarrow \infty$ in  \eqref{sumineq} yields
     \begin{align*}
    \sum_{k=0}^\infty\dfrac{\sigma_k}{2}||x_{k+1}-x_k||^2<+\infty
\end{align*}
since $\{\delta_k\}_k$ is summable. Therefore,  $||x_{k+1}-x_k|| \rightarrow 0 $
since  $\sigma_k$ is bounded from below by $1/\bar{\beta}-L>0$.

Now, we prove that  any cluster point $x_\infty$ is stationary for \eqref{sip}. Let $\{k_j\}_j$ be a subsequence of indices such that $x_{k_j}\rightarrow x_\infty$. Exploiting the same arguments of the  proof of Theorem \ref{ThmAlg1}, there exits a $H>0$ such that $$|\varphi_{Y_{i_{k_j},k_j}}(x_{k_{j+1}})- \varphi_{Y_{i_{k_j},k_j}}(x_{k_j})|\le H||x_{k_{j+1}}-x_{k_j}||$$
    holds for every $j$. Moreover, since $\varphi_{Y_{i_{k_j},k_j}}(x_{k_j})=\varphi(x_{k_j})$
    thanks to steps 6 and 10 of the algorithm, we obtain
\begin{align*}
    \varphi(x_{k_{j+1}})\ge\varphi_{Y_{i_{k_j},k_j}}(x_{k_{j+1}})&\ge \varphi_{Y_{i_{k_j},k_j}}(x_{k_j})-H||x_{k_{j+1}}-x_{k_j}||
    =\varphi(x_{k_j})-H||x_{k_{j+1}}-x_{k_j}||.
\end{align*}
 Taking the limit  as $j \rightarrow \infty$ in the above inequalities yields
$$\varphi_{\hat{Y}}(x_\infty)=\limsup_{j \rightarrow \infty} \varphi_{Y_{i_{k_j},k_j}}(x_{k_{j+1}})=\lim_{k_j \rightarrow \infty} \varphi_{Y_{i_{k_j},k_j}}(x_{k_{j+1}})=\varphi({x_\infty}),$$
 where $\hat{Y}$ is the Kuratowski superior limit of the sequence $\{Y_{i_{k_j},k_j}\}_j$, thanks to Proposition \ref{ConvergenceGeneralized Gradients}.
The  first-order optimality conditions of the proximal  linearized subproblem \eqref{proxlinprob} read:
$$\dfrac{x_{k_j}-x_{k_j+1}}{\gamma_{k_j}} \in \nabla f(x_{k_j})+ \rho_{k_j} \partial \varphi^+_{Y_{i_{k_j},k_j}}(x_{k_j+1})+\mathcal{N}_X(x_{k_j+1}).$$
Therefore,  $\varphi_{Y_{i_{k_j},k_j}}(x_{k_j+1})\le 0$ holds whenever $k_j \ge \bar{k}$. Thus, the same argument of     Proposition \ref{PenaltyFeasibility} shows that      there exists $\lambda_{k_j} \le \rho_{k_j}$ such that $$\dfrac{x_{k_j}-x_{k_j+1}}{\gamma_{k_j}} \in \nabla f(x_{k_j})+ \lambda_{k_j} \partial \varphi_{Y_{i_{k_j},k_j}}(x_{k_j+1})+\mathcal{N}_X(x_{k_j+1})$$
with
$$\lambda_{k_j}\varphi_{Y_{i_{k_j},k_j}}(x_{k_j+1})=0.$$
By taking the Kuratowski superior limit in the above inclusion, we obtain
    \begin{align*}
        0 &\in \Limsup_{j\rightarrow \infty}(\nabla f(x_{k_j})+\lambda_{k_j} \partial \varphi_{Y_{i_{k_j},k_j}}(x_{k_j+1})+ \mathcal{N}_X(x_{k_j+1}))\\
        & \subseteq \Limsup_{j\rightarrow \infty}(\nabla f(x_{k_j}))+\Limsup_{k_j\rightarrow \infty}(\lambda_{k_j}) \Limsup_{k_j\rightarrow \infty}(\partial \varphi_{Y_{i_{k_j},k_j}}(x_{k_j+1}))\\
        &+ \Limsup_{j\rightarrow \infty}(\mathcal{N}_X(x_{k_j+1}))\\
        &\subseteq
        \nabla f(x_{\infty})+\Lambda_\infty \partial \varphi_{\hat{Y}}(x_\infty)+ \mathcal{N}_X(x_\infty),
    \end{align*}
    where $\Lambda_\infty$ is the Kuratowski superior limit of the sequence $\{\lambda_{k_j}\}_j$. Without loss of generality, we can assume that $\lambda_{k_j}$ converges to some $\lambda_\infty \in \Lambda_\infty$ since the sequence  $\{\lambda_{k_j}\}_j$ is upper bounded by $\rho_{\bar{k}}$. Therefore, taking the limit as $j \rightarrow \infty$ in both  $$\lambda_{k_j}\varphi_{Y_{i_{k_j},k_j}}(x_{k_j+1})=0 \quad and \quad \varphi_{Y_{i_{k_j},k_j}}(x_{k_j+1})\le 0$$
    yields $\lambda_\infty\varphi_{\hat{Y}}(x_\infty)=0$ with  $\varphi_{\hat{Y}}(x_\infty)\le 0$. Hence, $x_\infty$ is stationary for  (VP$_{\hat{Y}}$).
   Thus,  $x_\infty$ is stationary for
    \eqref{sip}
by  exploiting the equality $\varphi_{\hat{Y}}(x_\infty)=\varphi(x_\infty)$ and Lemma \ref{GradInclusion}.
 \end{proof}

\textbf{Remark 1:}  Algorithms~2 and~3 at lines~$5$ and~$5'$ require  computing a minimizer to  a nonconvex optimization problem  to achieve convergence. Actually,  stationary points  are  enough  whenever their  value is
  no larger than the value at the current iterate, i.e.,  $f(x_k)+\rho_k\varphi^+(x_k)$ for Algorithm~2 and $\nabla f(x_k)^T x_k+\rho_k\varphi^+(x_k)$ for Algorithm~3.\\
  
  {  If $f$ is convex and $g(\cdot,y)$ is convex for every $y\in Y$, the objective functions of the  problems at  lines~5 and~5$'$ are strongly convex.
 Therefore, each problem has a unique
stationary point, which is also its global minimizer. In this convex setting, the boundedness of the sequence of penalty parameters  guarantees the convergence of the whole sequence of  iterates if at least one stationary point for \eqref{sip} exists.}
{\begin{corollary}
    Assume  (A4) for Algorithm 2 while (A5) for Algorithm 3. Moreover, assume also that $f$ is convex and $g(\cdot,y)$ is convex for every $y \in Y$, and that there exists a stationary point for \eqref{sip}.   If the sequence of penalty parameters
     is bounded, then  the sequence of iterates converges to a  stationary point  for \eqref{sip}.
\end{corollary}
\begin{proof}
     We only show the convergence for Algorithm 3, as the proof can be easily adjusted for Algorithm 2. Consider the  set $Y_{i_k,k}$ obtained  at the last iteration $i_k$ of the inner loop of the outer iteration $k$.
    Notice that  $x_{k+1}=z_{i_k,k}$ holds 
    thanks to step 9 of the algorithm, hence $x_{k+1}$ is a minimizer for the problem
     \begin{align}
         \min_{x \in X} \nabla f(x_{k})^Tx+ \rho_k \varphi^+_{Y_{i_k,k}}(x)+\dfrac{1}{2\gamma_k}||x-x_k||^2.\label{proxlinprob}
     \end{align}
     Now, let     $(x^*,\lambda^*)$ be any pair where $x^*\in X $ is stationary point  for \eqref{sip} along with multiplier $\lambda^*\ge0$.
     Then, 
     the following chain  holds
     \begin{align}
     \nabla f(x_{k})^Tx_{k+1}+&\dfrac{1}{2\gamma_k}||x_{k+1}-x_k||^2
\le \notag\\
&\le \rho_k \varphi^+_{Y_{i_k,k}}(x_{k+1})+\nabla f(x_{k})^Tx_{k+1}+\dfrac{1}{2\gamma_k}||x_{k+1}-x_k||^2\notag\\
     & \le \rho_k \varphi^+_{Y_{i_k,k}}(x^*)+\nabla f(x_{k})^Tx^*+\dfrac{1}{2\gamma_k}||x^*-x_k||^2-\dfrac{1}{2\gamma_k}||x^*-x_{k+1}||^2\notag\\&
     =\nabla f(x_{k})^Tx^*+\dfrac{1}{2\gamma_k}||x^*-x_k||^2-\dfrac{1}{2\gamma_k}||x^*-x_{k+1}||^2,\notag
     \end{align}
     where the  first inequality follows from $\varphi^+_{Y_{i_k,k}}\ge0 $, the second  inequality  from the fact that $x_{k+1}$ is a minimizer for problem \eqref{proxlinprob} that has a strongly convex objective function of parameter $1/\gamma_k$  and the  equality  from $ \varphi_{Y_{i_k,k}}^+(x^*)=0,$ which holds  since $x^*$ is feasible for \eqref{sip} and therefore for (VP$_{Y_{{i_k},k}}$).
     Thus, we obtain
     \begin{align}\label{proxineq1}
     \nabla f(x_{k})^T(x_{k+1}-x^*) \le   \dfrac{1}{2\gamma_k}||x_{k}-x^*||^2-\dfrac{1}{2\gamma_k}||x_{k+1}-x_k||^2
         -\dfrac{1}{2\gamma_k}||x_{k+1}-x^*||^2  
     \end{align}
     by rearranging the first and the last terms of the above chain.
      Assumption (A5) guarantees that the inequality 
$$f(x_{k+1})-f(x_k)+\nabla f(x_{k})^T(x_k-x_{k+1})\le  \dfrac{L}{2}||x_{k+1}-x_k||^2 $$
holds by \cite[Lemma 1.2.3]{Nesterov2004}, while the  convexity of $f$ guarantees the inequality
$$f(x_k)-f(x^*)+\nabla f(x_{k})^T(x^*-x_{k})\le 0, $$
and summing them  yields
$$f(x_{k+1})-f(x^*)+\nabla f(x_{k})^T(x^*-x_{k+1})\le \dfrac{L}{2}||x_{k+1}-x_k||^2.$$
By summing  the above inequality with \eqref{proxineq1}, we obtain
    \begin{align*}
    f(x_{k+1})-f(x^*)\le  \notag \dfrac{1}{2\gamma_k}||x_{k}-x^*||^2-\dfrac{\sigma_k}{2}||x_{k+1}-x_k||^2
         -\dfrac{1}{2\gamma_k}||x_{k+1}-x^*||^2\
\end{align*}
where $\sigma_k=1/\gamma_k-L>0$, which yields 
\begin{align}\label{proxineq2}
    f(x_{k+1})-f(x^*)
    \le &  \dfrac{1}{2\gamma_k}||x_{k}-x^*||^2-\dfrac{1}{2\gamma_k}||x_{k+1}-x^*||^2.
\end{align} 
By exploiting the stationarity conditions for \eqref{sip},  we obtain
$$0 \in \partial f(x^*)+\lambda^* \partial \varphi(x^*)+N_X(x^*)=\partial (f+\lambda^* \varphi)(x^*)+N_X(x^*)$$
where the equality follows from the  Clarke regularity of $f$ and $\varphi$ (see \cite[Corollary 2]{clarke_optimization_1990}). Since $f+\lambda^*\varphi$ is convex, the above inclusion yields that $x^*$ mimimizes $f+\lambda^*\varphi$ over the convex set $X$. This, together with the complementarity condition $\lambda^*\varphi(x^*)=0$, yields  
$$f(x^*)= f(x^*)+\lambda^*\varphi(x^*) \le f(x_{k+1})+\lambda^*\varphi(x_{k+1}).$$
 Notice that
$$\varphi(x_{k+1})\le \varphi^+(x_{k+1})\le \varphi^+_{Y_{i_k,k}}(x_{k+1})+\delta_k$$
holds, where the second inequality is due to  step 8 of the algorithm.
 Hence, combining the  above  two chains yields
$$f(x^*)-f(x_{k+1})\le \lambda^*\varphi^+_{Y_{i_k,k}}(x_{k+1})+\lambda^*\delta_k.$$ 
 Since the sequence $\{\rho_k\}_k$ is bounded, there exists $\bar{k}$ such that $g(x_{k+1},w_{k+1})\le 0$ holds for every $k \ge \bar{k}$ since the penalty parameter is increased only when  $g(x_{k+1},w_{k+1})>0$ (step 16 of the algorithm). Hence, the above inequality becomes  $$f(x^*)-f(x_{k+1})\le \lambda^*\delta_k$$
  whenever $k \ge \bar{k}$, since  $\varphi_{Y_{i_k,k}}(x_{k+1})=g(x_{k+1},w_{k+1})$ holds by step 7 of the algorithm.
Therefore, summing the above inequality with \eqref{proxineq2} yields that  \begin{align} \label{proxineq3}
    0
    \le &  ||x_{k}-x^*||^2-||x_{k+1}-x^*||^2+2\lambda^*\gamma_k\delta_k
\end{align}
holds for every $k \ge  \bar{k}$. The sequence $\{\lambda^*\gamma_k\delta_k\}_k$ is summable since each term  is nonnegative and bounded from above by $\lambda^*\bar{\beta}\delta_k$. Thus, \cite[Proposition 3.3]{Combettes2001} guarantees  that the sequence $\{||x_k-x^*||\}_k$ has a finite limit. Hence, the sequence $\{x_k\}_k$ is bounded and  admits a cluster point $x_\infty$  that  is also stationary for \eqref{sip} along with a suitable  multiplier $\lambda_\infty\ge 0$ by Theorem 11. Since \eqref{proxineq3} holds for every stationary pair $(x^*,\lambda^*)$, it also holds  for $(x_\infty,\lambda_\infty)$.
  Hence, the sequence $\{||x_k-x_\infty||\}_{k}$ goes to zero, since it has a finite limit by \cite[Proposition 3.3]{Combettes2001} and $x_\infty$ is a cluster point for $\{x_k\}_k$.
\end{proof}}

\textbf{Remark 2:}
Both the proximal algorithms can be modified by allowing  indices to be discarded  from the index sets $Y_k$.  Following the  strategy introduced in~\cite{blankenship_infinitely_1976}, it is not necessary to retain all the  indices previously generated since one may discard the constraints that are inactive at the current iterate. Accordingly, line~10 in Algorithm 2 or 3
may be replaced by
\begin{align*}
10':\qquad
Y_{k+1}
&:=
Y_0
\cup
\{ \
y\in Y_k
\ : \
g(x_{k+1},y)\geq 0
\ \}
\cup
\{\
w_{k+1},y_{k+1}
\ \}.
\end{align*}
By removing inactive  constraints the size of the resulting subproblems is controlled  without affecting convergence. In fact,  the same proofs {of Theorem 11 and Corollary 12} still work by taking
\begin{align*}
W_k
&:=
Y_0
\cup
\{\
y\in Y_k
\ :\
g(x_{k+1},y)\geq 0
\ \}
\cup
\{\
w_{k+1}
\ \},
\end{align*}
with no other modification.

The convergence results for the   algorithms rely on the boundedness of the penalty parameters. When the penalty parameter grows to infinity,   the methods either may not  converge or    get stuck at   infeasible points that violate  eNNAMCQ. Anyway, this can not happen  when the feasible region  $X$ is compact and every point satisfies eNNAMCQ since      the penalty parameter remains bounded.

\begin{theorem}
     If $X$ is compact and eNNAMCQ holds at every point of $X$, then the sequences of penalty parameters generated by  Algorithms 1 and 2 are bounded. If moreover (A5) holds, the boundedness of the sequence of penalty parameters is guaranteed also   for Algorithm 3.
\end{theorem}
\begin{proof}
First, we prove the theorem for Algorithm \ref{alg1}. Ab absurdo, suppose that the algorithm  produces an unbounded  sequence of penalty parameters. Then,  the compactness of $X$ guarantees that  there exists a subsequence of indices $\{k_j\}_j$ such that $x_{k_j}\rightarrow x_\infty,$ $\rho_{k_j}\rightarrow +\infty$ and $\varphi_{Y_{k_j}}(x_{k_j})>0$ with $\varphi_{Y_{k_j}}(x_{k_j}) \rightarrow\varphi(x_\infty)$. Moreover,
    $$0\in \partial f(x_{k_j})+\rho_{k_j}\partial\varphi_{Y_{k_j}}(x_{k_j})+\mathcal{N}_X(x_{k_j})$$
    holds since each point is stationary for the penalized problem (VP$_{Y_{k_j}}^{\rho_{k_j}}$).
    Then, dividing by $\rho_{k_j}$ and taking the Kuratowski superior limit yields
    \begin{align*}
    0 & \in  \Limsup_{{k_j}\rightarrow \infty}(\partial f(x_{k_j})/\rho_{k_j}  +  \partial \varphi_{Y_{{k_j}}}(x_{{k_j}}) + \mathcal{N}_X(x_{{k_j}})) \\
   & \subseteq  \Limsup_{{k_j}\rightarrow \infty}(\partial f(x_{k_j})/\rho_{k_j} )+\Limsup_{{k_j}\rightarrow \infty}(\partial \varphi_{Y_{k_j}}(x_{{k_j}})) + \Limsup_{{k_j}\rightarrow \infty}(\mathcal{N}_X(x_{{k_j}})) \\
   &\subseteq \partial \varphi_{\hat{Y}}(x_\infty) + \mathcal{N}_X(x_\infty)
    \subseteq \partial \varphi(x_\infty) + \mathcal{N}_X(x_\infty),
\end{align*}
which implies that $x_\infty$ violates eNNAMCQ, which gives the contradiction.

Now,  consider Algorithm 2. First of all, note that the assumption $(A4)$ is guaranteed since $X$ is compact. Ab absurdo, suppose that the algorithm  produces an unbounded  sequence of penalty parameters. Then,  the compactness of $X$ guarantees that  there exists a subsequence of indices $\{k_j\}_j$ such that $x_{k_j+1}\rightarrow x_\infty$ , $\rho_{k_j}\rightarrow +\infty$ and $\varphi_{Y_{k_j}}(x_{k_j+1})>0$ with $\varphi_{Y_{k_j}}(x_{k_j+1}) \rightarrow\varphi(x_\infty)$. Moreover,
    $$\dfrac{x_{k_j}-x_{k_j+1}}{\gamma_{k_j}}- \partial f(x_{k_j+1})\in \rho_{k_j}\partial\varphi^+_{Y_{k_j}}(x_{k_j+1})+\mathcal{N}_X(x_{k_j+1})$$
     holds since each point is an optimal solution to the proximal subproblem at line 5 of the algorithm. Then,  dividing by $\rho_{k_j}$ and taking the Kuratowski superior  limit as $j \rightarrow \infty$ yields the contradiction as in the previous case since  $\{x_{k_j}-x_{k_j+1}\}_j$ is bounded due the compactness of $X$. The same proof works  for Algorithm 3  with  $\nabla f(x_{k_j+1})$ in place of $\partial f(x_{k_j+1})$ under assumption (A5).
\end{proof}

\section{Numerical tests}

 The algorithms of the previous section are  compared with recent  methods for semi-infinite optimization that share similar features.  For this reason, the numerical tests are organized into three sets of experiments. First, we compare the basic  penalty and proximal algorithms with  their variants in order to determine the method that performs the best.  This happens to be Algorithm 3 with the constraint-dropping strategy. Therefore, this latter algorithm is  then  compared  with the integral penalty method~\cite{xu_solving_2014} and  with  other cutting surface algorithms~ \cite{blankenship_infinitely_1976,mehrotra_cutting_2014,auslender_comparative_2015} for convex semi-infinite problems.

All the experiments were performed on a laptop equipped with an Intel Core i7-8750H processor running at $2.20$ GHz and a $64$-bit Windows operating system. The algorithms were implemented in Python~3.8 using CVXPY as the modeling framework. Convex subproblems were solved with SCS~\cite{SCS}, whereas nonconvex subproblems were handled by the DIRECT implementation available in SciPy~\cite{DIRECT,Scipy}.

\subsection{Comparison between  the proposed algorithms}
Five algorithms have been
considered: \texttt{Penalty}, \texttt{Prox}, \texttt{ProxLin}, \texttt{ProxDrop}, and \texttt{ProxLinDrop}. \texttt{Penalty} is  Algorithm 1 described in Section 4, whereas  \texttt{Prox} and  \texttt{ProxLin} are respectively Algorithm 2 and 3. Lastly, the suffix  \texttt{Drop} indicates the variant of both the latter algorithms where the indices corresponding to the  inactive constraints are discarded.

The algorithms have been  tested on Problem~1, a minimax problem that generalizes the test problem considered in~\cite{mehrotra_cutting_2014} by considering arbitrary coefficients in the  semi-infinite constraint. This  is a convex semi-infinite program and admits a Slater point. Hence, eNNAMCQ holds globally, and the penalty parameters generated by the algorithms remain bounded. As a consequence, the convergence results of the  previous section  hold.

\medskip

\textbf{Problem 1} (Chebyshev-type approximation):
\begin{align*}
\min \quad
& \eta\\
\textnormal{s.t.}\quad
& x\in[-1,1]^n,
& \sum_{i=1}^n
\left(
\alpha_i x_i-\frac{\beta_i}{n}
-\sin(2\pi t+\gamma_i)
\right)^2
\leq \eta,
\qquad
\forall t\in[0,1].
\end{align*}

\medskip

The tests were run on  $100$ instances.
The dimension was set to $n=50$.  The parameters $\alpha_i$, $\beta_i$, and $\gamma_i$ were sampled independently from the uniform distribution on $[i,i+1]$ for  every $i$. The initial penalty parameter and the update factor were set to $\rho_0=1$ and $\sigma=2$ for all the algorithms. Finally,  $\gamma_k=1$ was chosen as a fixed value and the tolerance sequence was defined by $\delta_k=0.9^k$
for the proximal algorithms.

The stopping criterion for the \texttt{Penalty} algorithm was chosen as $\varphi(x_k)\leq 10^{-3}$, whereas the  stopping criterion for the proximal methods additionally required $||x_k-x_{k-1}||\leq 10^{-3}$.

Table~\ref{tab:mean_prox} reports the means of the computational time, a measure of   feasibility, the number $k$ of  iterations, the number of evaluations of the  value function $\varphi$, and the number of penalty updates. The feasibility was measured by evaluating the value function $\varphi$ at the last iterate $x_k$. Indeed,  this quantity is positive when the  solution violates the semi-infinite constraint of \eqref{sip}, while a negative value implies that the solution is feasible. Notice that, the numbers of iterations and evaluations  is different in the proximal methods since  each  iteration may require multiple evaluations of the value function due to the inner loop.

\begin{table}[h]
\centering
\begin{tabular}{lrrrrr}
\hline
algorithm & time & feasibility & iter & eval & $\rho$ incr \\
\hline
Penalty     & 16.59  & $5.95 \times 10^{-4}$ & 11.73 & 11.73 & 10.69 \\
Prox        & 17.13  & $7.22 \times 10^{-4}$ & 12.57 & 16.11 & 9.32  \\
ProxLin     & 15.69  & $4.19 \times 10^{-4}$ & 8.43  & 13.02 & 5.49  \\
ProxDrop    & 12.14  & $5.16 \times 10^{-4}$ & 18.50 & 22.38 & 9.76  \\
ProxLinDrop & 10.91  & $4.70 \times 10^{-4}$ & 11.19 & 15.90 & 6.37  \\
\hline
\end{tabular}
\caption{Mean results for the proposed algorithms.}
\label{tab:mean_prox}
\end{table}

The results show that the linearization of the objective  function reduces both the number of iterations and penalty updates. Removing inactive constraints mainly reduces computational time, although \texttt{ProxDrop} requires more iterations and value-function evaluations. Combining the two features gives the best overall performance: \texttt{ProxLinDrop} has the lowest mean computational time and the  constraint violation
remains comparable with that of the other variants.

\begin{figure}[htbp]
\centering
\includegraphics[width=0.7\linewidth]{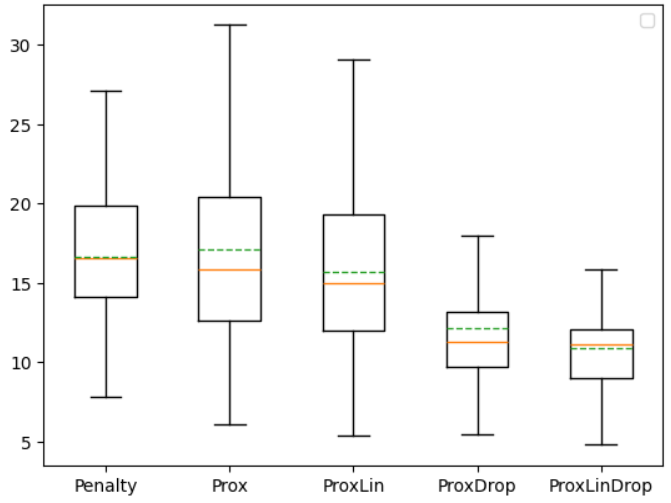}
\caption{Performance comparison of the algorithms on Problem~1.}
\label{fig1}
\end{figure}

Figure~\ref{fig1} reports the  box plot of the  computational times. The orange line represents the median, whereas the dashed green line represents the mean. The spreads of the results  of \texttt{Penalty}, \texttt{Prox}, and \texttt{ProxLin} are wider than those of the index-dropping variants. In particular, \texttt{ProxLinDrop} achieves both the lowest mean and median computational times and a relatively narrow distribution.  It is therefore selected for the subsequent comparisons.

\subsection{Comparison with the integral penalty method}

We next compared \texttt{ProxLinDrop} with the integral penalty method proposed in~\cite{xu_solving_2014}, referred to as \texttt{SPG} (Smoothed Projected Gradient).  This  method replaces the semi-infinite constraint with an integral-type constraint and applies a penalty approach to the resulting problem. The comparison relies on  Problems~2--4 taken from \cite{xu_solving_2014}. Unlike Problem~1, these problems involve only two  variables but the objective function and/or the functions $g(\cdot,y)$ are not convex.

\textbf{Problem 2}
(sum-of-square objective):
\begin{align*}
    \min \quad & (x_1-2x_2+5x_2^2-x_2^3-13)^2 + (x_1-14x_2+x_2^2+x_2^3-29)^2\\
    \textnormal{s.t.}\quad & x\in \mathbb{R}^2, x_1^2+2x_2y^2+\exp(x_1+x_2)-\exp(y)\le 0, \ \forall y\in[0,1]
\end{align*}

\textbf{Problem 3}
(polynomial functions):
\begin{align*}
    \min \quad & 1/3x_1^2 +1/2x_1 + x_2^2\\
    \textnormal{s.t.}\quad  & x\in [-10,10]^2,\  (1-x_1^2y^2)^2-x_1y^2-x_2^2+x_2\le 0,\ \forall y \in [0,1]
\end{align*}

\textbf{Problem 4} (exponential functions):
\begin{align*}
    \min \quad & 1.5^2\exp(x_1)+\exp(x_2)\\
     \textnormal{s.t.}\quad &   x\in [-10,10]^2,\  y-\exp(x_1+x_2)\le 0, \ \forall y \in [0,1]
\end{align*}
\begin{table}[h]
\centering

\begin{tabular}{clrrr}
\hline
problem &algorithm & time & feasibility & optimal \\
\hline
\multirow{2}{*}{2} & ProxLinDrop     & 0.76  & $-3.14 \times 10^{-7}$ &  97.16 \\
                   & SPG & 26.59 & $-1.19 \times 10^{-2}$ & 97.16 \\
\hline
\multirow{2}{*}{3} & ProxLinDrop     & 2.38  & $-1.68 \times 10^{-8}$ & 0.19 \\
                   & SPG & 33.13 & $-8.98 \times 10^{-4}$ & 0.20 \\
\hline
\multirow{2}{*}{4} & ProxLinDrop     & 2.10  & $3.62 \times 10^{-5}$  & 3.00 \\
                   & SPG & 11.17 & $-1.12 \times 10^{-2}$ & 3.02 \\
\hline
\end{tabular}
\caption{comparison between ProxLinDrop and  SPG}
\label{tab:integral-penalty}
\end{table}

Table~\ref{tab:integral-penalty} reports the computational time, the measure of  feasibility, and the optimal value. \texttt{ProxLinDrop} requires less computational time on all three problems.
The two algorithms compute nearly identical optimal values, with the integral penalty method always returning feasible points  while \texttt{ProxLinDrop}
slightly violates some  constraints only in one problem. A substantial  advantage in the computational time is observed for \texttt{ProxLinDrop}, possibly due to the increasing cost of evaluating the integral as the smoothing parameter of \texttt{SPG} grows.

\subsection{Comparison with other cutting surface methods}

We finally compared \texttt{ProxLinDrop} with three algorithms for convex semi-infinite optimization: the classical cutting surface algorithm~\cite{blankenship_infinitely_1976}, denoted by \texttt{CS}; the central cutting surface algorithm~\cite{mehrotra_cutting_2014}, denoted by \texttt{CCS}, and  the Remez penalty and smoothing algorithm proposed in~\cite{auslender_comparative_2015}, denoted by \texttt{Remez}. The comparison relies on  Problem~1 and two additional families of convex test problems that are reported below.

Problem~5 is a linear semi-infinite problem that extends the test problem proposed in~\cite{WeiHaskellZhao2020} from fixed  to random coefficients. The dimension of the index set is twice that of the variables, making the problem computationally demanding even for small values of $n$.

\textbf{Problem 5} (linear semi-infinite problem):
\begin{align*}
&\min c^Tx\\
& s.t.\ x\in [-2,2]^n, \max_{i=1,\dots,2n}\max_{||y_i||\le 1}(A_i+\delta_iy_i)^Tx+b\le 0
\end{align*}

At last, we consider a new problem that has a linear objective and a quadratic constraint indexed by the $n$-dimensional unit ball. The indices $y$ appear in the quadratic,  linear and exponential terms, giving rise to an infinite  family of  convex nonlinear
constraints.

\textbf{Problem 6 }(quadratic semi-infinite constraint):
\begin{align*}
&\min c^Tx\\
&s.t. \   x \in [-2,2]^n, \   x^T(Q+\text{diag}(y))x + y^Tx -e^{q^Ty} \le 0, \ ||y||\le 1,
\end{align*}

The dimensions were set to $n=50$ for Problem~1 and 6 while to  $n=10$ for Problem~5, and $100$ instances were generated for each problem with parameters sampled independently from uniform distribution. In Problem~1 $\alpha_i$, $\beta_i$, and $\gamma_i$ were sampled  from  $[i,i+1]$ for each $i$. In Problem~5 we set
\[
A=
\begin{bmatrix}
-I\\
I
\end{bmatrix},
\]
where $I$ is the identity matrix of dimension $n$, while $\delta_i$, $c$, and $b$ were sampled from $[0.1,0.3]$, $[-5,5]^n$, and $[-2,2]^{2n}$, respectively. In Problem~6, $q$ and $c$ were sampled from $[0.2,0.6]^n$ and $[-5,5]^n$, respectively. Finally,  $Q\in\mathbb{R}^{n\times n}$ was randomly generated to satisfy $Q\succeq I$ by considering
$$Q=P^TP+I,$$
where $P$ was sampled   in $[-1,1]^{n\times n}$.

\begin{table}[htbp]
\centering
\begin{tabular}{lrrrr}
\hline
algorithm & time & feasibility & iter & eval \\
\hline
ProxLinDrop & 10.91  & $4.70 \times 10^{-4}$ & 11.19 & 15.90 \\
CCS        & 33.19 & $1.20 \times 10^{-3}$ & 13.25 & 13.25 \\
CS          & 39.03 & $8.63 \times 10^{-4}$ & 17.43 & 17.43 \\
Remez       & 30.71 & $8.35 \times 10^{-4}$ & 13.81 & 13.81 \\
\hline
\end{tabular}
\caption{Mean results on  Problem 1}
\label{tab:mean1}
\end{table}

\begin{table}[htbp]
\centering
\begin{tabular}{lrrrr}
\hline
algorithm & time & feasibility & iter & eval \\
\hline
ProxLinDrop & 47.37  & $8.54 \times 10^{-5}$ & 50.17 & 66.53 \\
CCS        & 72.79  & $8.27 \times 10^{-4}$ & 55.37 & 55.37 \\
CS          & 57.27  & $7.21 \times 10^{-4}$ & 48.10 & 48.10 \\
Remez       & 95.22  & $5.19 \times 10^{-4}$ & 73.77 & 73.77 \\
\hline
\end{tabular}
\caption{Mean results on Problem 5}
\label{tab:mean2}
\end{table}

\begin{table}[htbp]
\centering
\begin{tabular}{lrrrr}
\hline
algorithm & time & feasibility & iter & eval \\
\hline
ProxLinDrop & 6.32   & $1.93 \times 10^{-4}$ & 9.66  & 11.66 \\
CCS      & 109.69 & $7.51 \times 10^{-4}$ & 18.56 & 18.56 \\
CS         & 22.51  & $5.45 \times 10^{-4}$ & 7.86  & 7.86  \\
Remez       & 412.93 & $8.98 \times 10^{-4}$ & 38.80 & 38.80 \\
\hline
\end{tabular}
\caption{Mean results on Problem 6}
\label{tab:mean3}
\end{table}

Tables~\ref{tab:mean1}--\ref{tab:mean3} show that \texttt{ProxLinDrop} consistently achieves the lowest mean computational time and the smallest mean constraint violation across all the three problems.
Although the number of evaluations of  the value function $\varphi$ performed  by \texttt{ProxLinDrop}  is comparable with the other methods, it is never the lowest.
This suggests that the gap in the   computational time is mainly due to subproblems that are easier to solve. Indeed, the proximal regularization makes the subproblems strongly convex, while the index-dropping strategy limits the growth of the number of constraints.

\begin{figure}[htbp!]
    \centering
    \includegraphics[width=0.7\linewidth]{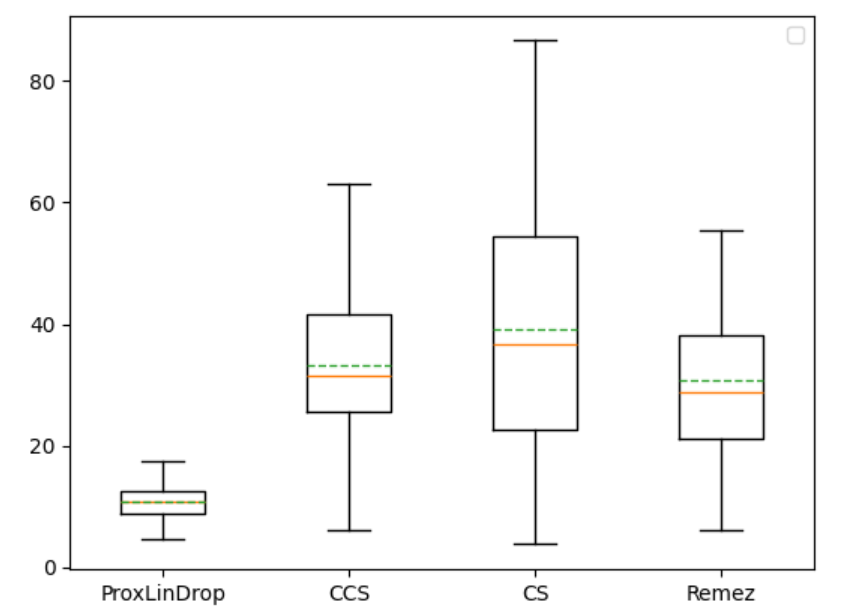}
    \caption{Box plot of the computational times  on Problem 1}
    \label{fig2}
\end{figure}

\begin{figure}[htbp!]
    \centering
    \includegraphics[width=0.7\linewidth]{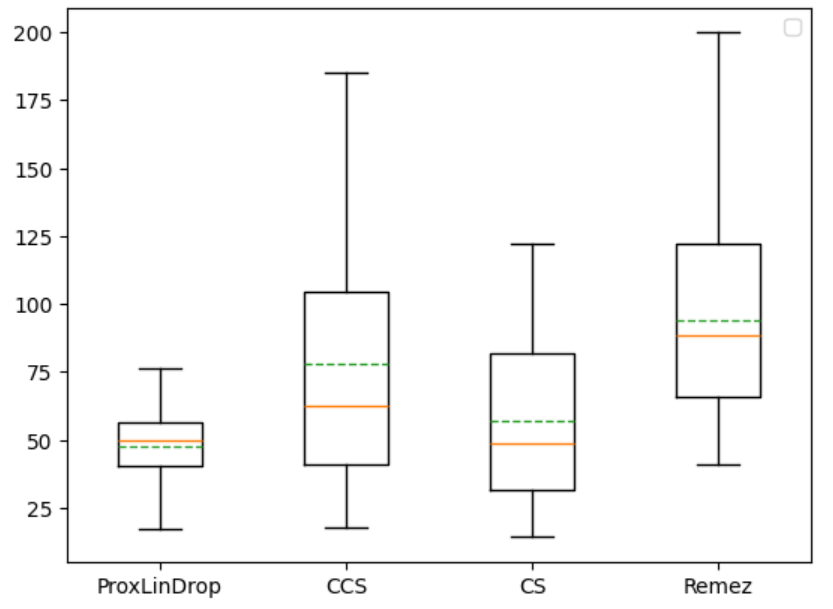}
    \caption{Box plot of the computational times  on Problem 5}
    \label{fig3}
\end{figure}

\begin{figure}
    \centering
    \includegraphics[width=0.7\linewidth]{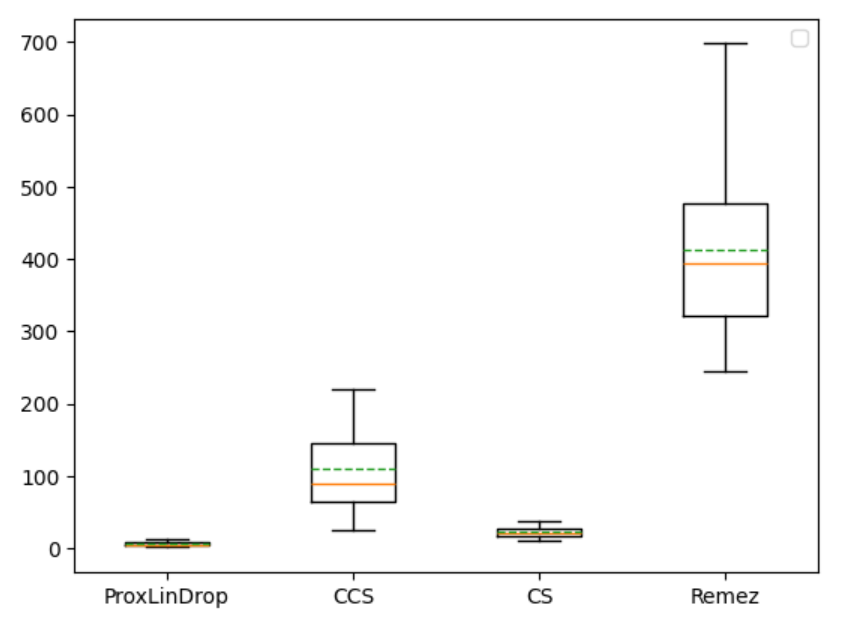}
    \caption{Box plot of the computational times  on Problem 6}
    \label{fig4}
\end{figure}

Figures~\ref{fig2}--\ref{fig4} report  the corresponding box plots  of the  computational time, which show that \texttt{ProxLinDrop} has also the lowest and most concentrated spread over all the tests. The box plots therefore  support the conclusions drawn from the mean results and indicate that \texttt{ProxLinDrop} consistently performs better than the other  methods throughout all the tests.

\section{Conclusions}

In this paper, algorithms were devised by combining cutting surface approximations with an adaptive penalty scheme and proximal regularization in order to  compute stationary points of a nonsmooth semi-infinite problem. The convergence of the methods was analyzed with respect to the behavior of the penalty parameters, and an exact penalty result was established under the compactness of  the feasible region \(X\)  and eNNAMCQ holding  at every point.
Numerical results show that Algorithm 3  combined with the index-dropping strategy   yields the best    performance among the proposed  algorithms and that it also  compares favorably with similar methods from the literature. Further work can be carried out in some directions.

When the sequence of penalty parameters is not bounded, the algorithms may get stuck at infeasible point that don't satisfy eNNAMCQ. One possibility to overcome this issue is to relax the constraints by allowing a small violation. This guarantees the existence of a Slater point, which is enough to restore exact penalization  in the convex setting.  However, this regularity is not sufficient to guarantee exact penalization in  nonconvex settings, and additional convexification techniques may be needed.

The approximate problems could also be solved inexactly by accepting  points that satisfy the stationarity conditions within a prescribed tolerance. This could reduce the computational effort as the number of considered  constraints  increases, but would require analyzing how the chosen tolerance affects the generation of new constraints and, consequently, convergence.

 The exploitation of smoothing techniques in combination with our approach could also be explored. Smooth approximations of the value function could be considered, while the positive-part penalty term could be replaced by a smooth penalty function. This would require analyzing the interaction between the smoothing and penalty parameters.

The framework could also be extended to generalized semi-infinite programming (GSIP). Although some  discretization and cutting surface methods have been devised for GSIP \cite{still_discretization_2001,Still2001}, the dependence of the index set on the variables makes their construction and convergence analysis considerably more complicated. At the same time, the study of exact penalization  for GSIP has been mainly   theoretical \cite{yang_optimality_2016}, hence it is not known how these techniques can be effectively combined with cutting surface  approximations. Extending the present framework to this more complex setting therefore requires further theoretical and algorithmic developments.

\backmatter

\bmhead{Acknowledgements}

The authors are members of the Gruppo Nazionale per l’Analisi
Matematica, la Probabilit\`a e le loro Applicazioni (GNAMPA - National Group for Mathematical
Analysis, Probability and their Applications) of the Istituto Nazionale di Alta
Matematica (INdAM - National Institute of Higher Mathematics).

\section*{Declarations}

\noindent\textbf{Data availability}

No data sets were generated during the current study. The used python codes are available on GitHub at \url{https://github.com/RIKTOMA/SIP}.

\bibliography{bibliography}
\end{document}